\documentclass[review]{elsarticle}
\usepackage[margin=2.5cm]{geometry}
\usepackage{hyperref}
\usepackage[]{algorithm2e}
\usepackage{enumitem}

\journal{  }

\usepackage{eurosym}
\usepackage{amsmath}
\usepackage{amsfonts}
\usepackage{array}
\usepackage{amssymb}
\usepackage{textcomp}
\usepackage{changepage}
\usepackage{graphicx}
\usepackage{caption}
\usepackage{subcaption}
\usepackage{natbib}
\usepackage{xcolor}
\usepackage{booktabs}
\usepackage{multirow}

\usepackage[authormarkup=none]{changes}
    \definechangesauthor[name=Sauleh,color=blue]{Sauleh}

\biboptions{authoryear}

\usepackage{etoolbox}
\makeatletter
\patchcmd{\ps@pprintTitle}
  {Preprint submitted to}
  {}
  {}{}
\makeatother

\begin{document}

\begin{frontmatter}

\title{ Modelling an electricity market oligopoly with a competitive fringe and generation investments}


\author[mymainaddress]{Mel T. Devine\textsuperscript{a,}\corref{mycorrespondingauthor}}
\cortext[mycorrespondingauthor]{Corresponding author}
\ead{mel.devine@ucd.ie}

\author[ucdaddress_d]{Sauleh Siddiqui\textsuperscript{c,}}

\address[ucdaddress_b]{College of Business, University College Dublin, Ireland}

\address[mymainaddress]{School of Electrical and Electronic Engineering, University College Dublin, Ireland}

\address[ucdaddress_c]{Department of Civil and Systems Engineering, Johns Hopkins University, Baltimore, USA}

\address[ucdaddress_d]{German Institute for Economic Research (DIW Berlin), Mohrenstr. 58, 10117, Berlin, Germany}

\begin{abstract}
Market power behaviour often occurs in modern wholesale electricity markets. Mixed Complementarity Problems (MCPs) have been typically used for computational modelling of market power when it is characterised by an oligopoly with competitive fringe. However, such models can lead to myopic and contradictory behaviour. Previous works in the literature have suggested using conjectural variations to overcome this modelling issue. We first show however, that an oligopoly with competitive fringe where all firms have investment decisions, will also lead to myopic and contradictory behaviour when modelled using conjectural variations. Consequently, we develop an Equilibrium Problem with Equilibrium Constraints (EPEC) to model such an electricity market structure. The EPEC models two types of players: price-making firms, who have market power, and price-taking firms, who do not. In addition to generation decisions, all firms have endogenous  investment decisions for multiple new generating technologies. The results indicate that, when modelling an oligopoly with a competitive fringe and generation investment decisions, an EPEC model can represent a more realistic market structure and overcome the myopic behaviour observed in MCPs. The EPEC considered found multiple equilibria for investment decisions and firms’ profits. However, market prices and consumer costs were found to remain relatively constant across the equilibria. In addition, the model shows how it may be optimal for price-making firms to occasionally sell some of their electricity below marginal cost in order to de-incentivize price-taking firms from investing further into the market.  Such strategic behaviour would not be captured by MCP or cost-minimisation models.
\end{abstract}

 
\end{frontmatter}

\section{Introduction} \label{sec:intro}
Electricity market modelling is an area of research that has attracted much attention in the Operations Research literature. Optimisation and equilibrium models in particular have been extensively used to better understand the behaviour of electricity generators. Such tools provide insights from planning, operations and regulatory perspectives. Regulators may use them to monitor market inefficiencies, profit-maximising generators may use them to gain insights on possible trading strategies while policy-makers may use them to test the impact of different proposed policy mechanisms.

Since the 1980s, countries have been deregulating their electricity markets with the intention of splitting ownership of market activities \citep{pozo2017basic}. Governments' goals are to foster competition, increase market efficiencies and thus reduce consumer costs. As a result, individual market participants, also known as market players, have been behaving selfishly by seeking to independently maximise their profits \citep{facchinei2007finite}.

Deregulation has resulted in many electricity markets showing evidence of market power \citep{lee2016nash}. Market power is present when one (or more) seller(s) in the market can strategically maximise their profits by influencing the selling price through the quantity they supply to the market. When such behaviour is not present in the market, the market is perfectly competitive. Accurately modelling market power in electricity markets is a challenging area of operations research. However, there has been many electricity market models that have incorporated market power. For a comprehensive review of electricity market models that incorporate market power, we refer the reader to \cite{pozo2017basic}. More recent examples in the operations research literature include  \cite{devine2018examining} who use a Mixed Complementarity Problem (MCP) to study different consumer led load shedding strategies. MCPs solve multiple constrained optimisation problems simultaneously and in equilibrium and they allow players with market power to be modelled as Cournot players. 

\cite{fanzeres2019robust} proposed a Mathematical Program with Equilibrium Constraints (MPEC) for strategic bidding in an electricity market. An MPEC model solves a bi-level optimisation problem which is a mathematical program where one or more optimisation problems are embedded within another optimisation problem. The outer optimisation is the upper-level optimisation while the inner optimisations, which are represented in the outer problem as constraints, are the lower-level optimisation problems. In an electricity market setting, MPEC problems can be used to model markets where a single player has market power. The upper level represents the optimisation problem of the player who has market power while the lower level problems represents the problems of players that do not have market power.

\cite{steeger2017dynamic} present a methodology that combines Lagrangian relaxation and nested Benders decomposition to model a single hydro producer with market power. Similarly, \cite{habibian2019multistage} us an optimisation-driven heuristic approach to model a large electricity consumer with market power. In both of these works, only one market participant has a strategic advantage.

The papers in the previous paragraphs do not consider strategic behaviour when the overall market was characterised by an by oligopoly with competitive fringe. Such a market structure occurs when more than one generator (the oligopolists) have market power and at least one generator does not (the competitive fringe). Many modern electricity markets are characterised by an oligopoly with competitive fringe \citep{bushnell2008vertical, walsh2016strategic}. However, when it comes to the energy market modelling literature, such markets structures are under-represented. 

Some exceptions include \cite{huppmann2013endogenous} and \cite{ansari2017opec}, who use a Mixed Complemenatrity Problem (MCP) to model an oligopoly with competitive fringe in international oil market contexts. Similarly, \cite{devine2019role} develop a MCP model of an oligopoly with competitive fringe to investigate the impact demand response has on market power in an electricity market. \cite{huppmann2013endogenous} highlights how modelling an oligopoly with competitive fringe using a MCP can lead to myopic, counter-intuitive and thus unrealistic optimal decisions from the oligopolists. In a MCP framework, each oligopolist optimises its own position but does not take into account the optimal reaction of the competitive fringe. The oligopolists reduce their generation levels with the intention of increasing the market price and hence increasing their overall profits. However, when the oligopolists do this, the competitive fringe increase their generation and fill the generation gap. This results in market prices not increasing as the oligopolists would anticipate. \cite{huppmann2013endogenous} proposes using conjectural variations to overcome the issue. Conjectural variations make assumptions about how players react to other players' quantity changes and have been widely used in the energy market modelling literature \citep{egging2008complementarity, haftendorn2010modeling, huppmann2014market, baltensperger2016multiplicity, egging2016risks}. These assumptions allow the oligopolists to somewhat incorporate the reactions of the competitive into their decision making process. However, the resulting decisions from the oligopolists may not be necessarily be optimal. 

Neither \cite{huppmann2013endogenous}, \cite{ansari2017opec} nor \cite{devine2019role} consider investment in new generation decisions. Consequently, in this work, we show that when investment decisions are incorporated into a MCP model of an oligopoly with competitive fringe, conjectural variations still lead to contradictory behaviour from the oligopolists. Moreover, to overcome the modelling issue, we develop an Equilibrium Problem with Equilibrium Constraints (EPEC) model of an oligopoly with competitive fringe where both oligopolists and the competitive fringe have investment decisions. Each firm may initially hold, and invest in, multiple generating technologies.

An EPEC model solves multiple interconnected MPEC problems in equilibrium \citep{gabriel2012complementarity}. In this work, each MPEC represents the optimisation problem of a different electricity generating firm who has market power (also known as a price-making firm). The price-making firms each seek to maximise profits subject to capacity constraints. In addition, the equilibrium conditions representing the optimisation problems of the competitive fringe are embedded into each price-making firm's problem as constraints. Consequently, an EPEC approach can overcome the limiting assumptions associated with conjectural variations. Instead of making assumptions of how players react to other players' quantity changes, an EPEC model allows oligopolists to explicitly account for optimal reactions of the competitive fringe and thus to make optimal decisions by anticipating those reactions.

There has been many examples in the literature where EPECs have been used to model electricity markets. Many of the first EPEC models for electricity markets consider electricity generators in the upper level and an Independent System Operator (ISO) in the lower level \citep{hu2007using, ruiz2011equilibria, pozo2011finding}. Building on these works, \cite{wogrin2012capacity} and \cite{wogrin2013generation} develop an EPEC model that incorporates both capacity expansion and generation decisions amongst electricity generators. In the upper level, the generators decide investment decisions whilst accounting for operational decisions in the lower level. \cite{pozo2013if} and \cite{pozo2012three} developed a model similar to \cite{wogrin2012capacity} but add an extra level to the model; an ISO who makes transmission expansion decisions whilst accounting for generators' capacity investment and operational decisions. \cite{jin2013tri}  consider a similar EPEC model as well but, in contrast to \cite{pozo2013if} and \cite{pozo2012three}, model price-responsive demand and strategic interactions amongst the generators. 

\cite{kazempour2013generation} and \cite{kazempour2013equilibria} use EPEC models where the upper-level problem determines the optimal investment for strategic producers while lower-level problems represent different market clearing scenarios. Similarly, \cite{ye2017investigating} use a EPEC model to investigate the impact consumer led demand shifting has on market power and find that demand response can reduce the negative impacts of market power. The upper level again represents the producers' problems while the lower-level represent the market clearing process, in addition to the consumers decisions. An EPEC model is used in \cite{huppmann2015national} as part of a three-stage equilibrium model between a supra-national planner, zonal planners, and an ISO. \cite{moiseeva2017generation} develop an EPEC model that considers generators' operational decisions in the lower level and their ramping decisions in the upper level. More recently, \cite{guo2019electricity} introduce another EPEC model where the upper level maximises generators' decisions while the lower level represents an ISO. Interestingly, \cite{guo2019electricity} account for risk-averse decision making by incorporating Conditional Vale at Risk (CVaR) into their model.

Despite the rich literature of EPEC models for electricity markets, none of the aforementioned EPEC problems model a market characterised by an oligopoly with competitive fringe, where some generators have market power and other do not. The closest work to the present paper is \cite{zerrahn2017network}. They propose a three-stage game to model transmission network expansion in an imperfectly competitive market where some generators have market power while do not. They solve the model using backward induction. The third stage represents the problem of the ISO and the competitive fringe. The second stage represents the firms who have market power and thus account for the third stage. In the first stage, social welfare is maximised using network expansion decisions whilst accounting for the second and third stages. Significantly, we advance the work of \cite{zerrahn2017network} by including generation expansion/investments for generating firms in our model. 

We apply the EPEC developed in this work to an electricity market representative of the Irish power system in 2025 using data from \cite{lynch2019role} and \cite{bertsch2018analysing}. EPEC models can be challenging to solve computationally. We utilise the Gauss Seidel algorithm in to order use diagonalization for solving the EPEC, and we solve each individual MPEC using disjunctive constraints \citep{fortuny1981representation}. In addition, to improve computational efficiency, we utilise the approach developed in \cite{leyffer2010solving} to provide an initial strong stationary point of the EPEC to use as a starting point for our diagonlization algorithm. We solve the model numerically as it is too large to be solved in closed form, which is another contribution of this work. A closed-form solution is possible using standard techniques but we combine two techniques from the literature in order to solve our problem. 

Our results show that may it be optimal for generating firms with market power to occasionally operate some of their generating units at a loss. The driving factor behind this model outcome is the fact we allow both price-making and price-taking firms to make investment decisions. Price-taking firms' ability to invest further into the market motivates the price-making firms to depress prices in some timepoints. This reduces the revenues price-taking firms could make from new investments and thus prevents them from making such investments. Such strategic behaviour would not be captured by MCP or cost-minimisation models. Consequently, this result highlights the suitability of the EPEC modelling approach and the importance of including investment decisions in models of oligopolies with competitive fringes.

The remainder of this paper is structured as follows: firstly, in Section \ref{sec:input_data}, we describe the model data inputs. Secondly, in Section \ref{sec:results_mcp} we demonstrate the naivety of using a MCP to model an oligopoly with a competitive fringe where both price-making and price-taking firms have investment decisions. Thirdly, in Section \ref{sec:model}, we introduce the EPEC model. In Sections \ref{sec:dis} and \ref{sec:con}, we provide some discussion and conclusions, respectively. Finally, in \ref{sec:app_pm}, we provide additional material related to the case study.

\section{Input data}\label{sec:input_data}

\begin{table}[htp!]
\footnotesize
\caption{Indices and sets. }\label{tab:sets}
\begin{tabular}{rl}
\hline
$f\in F$ &Generating firms \\
$t \in T$& Generating technologies\\
$p \in P$& Time periods\\
 \hline
\end{tabular}%
\end{table}

\begin{table}[htp!]
\footnotesize
\caption{Variables.}\label{tab:variables}
\begin{tabular}{rl}
\multicolumn{2}{l}{Price-taking firms' primal variables}\\
\hline
$gen^{\text{PT}}_{f,t,p}$& Forward generation from price-taking firm $f$ with technology $t$ in period $p$\\
$inv^{\text{PT}}_{f,t}$&Investment in new generation capacity (technology $t$) for price-taking firm $f$\\
\multicolumn{2}{l}{Price-making firms' primal variables}\\
\hline
$gen^{\text{PM}}_{l,t,p}$& Forward generation from price-making firm $l$ with technology $t$ in period $p$\\
$inv^{\text{PM}}_{l,t}$&Investment in new generation capacity (technology $t$) for price-making firm $l$\\
\multicolumn{2}{l}{Dual variables}\\
\hline
$\gamma_{p}$ & System price for time period $p$\\
$\lambda^{\text{PT}}_{f,t,p}$&Lagrange multiplier associated with price-taking firm $f$'s capacity constraint for technology $t$ and timestep $p$\\
$\lambda^{\text{PM}}_{l,t,p}$&Lagrange multiplier associated with price-making firm $l$'s capacity constraint for technology $t$ and timestep $p$\\
 \hline
\end{tabular}%
\end{table}

\begin{table}[htp!]
\footnotesize
\caption{Parameters.}\label{tab:parameters}
\begin{tabular}{cp{12.5cm}}
\hline
$CAP^{\text{PT}}_{f,t}$& Initial generating capacity for price-taking firm $f$ with technology $t$\\
$CAP^{\text{PM}}_{f,t}$& Initial generating capacity for price-taking firm $l$ with technology $t$\\
$A_{p}$ & Demand curve intercept for timestep $p$\\
$B$ & Demand curve slope\\
$C^{\text{GEN}}_{t}$ & Marginal generation cost for technology $t$ \\
$IC^{\text{GEN}}_{t}$& Investment in generating technology $t$ cost  \\
$W_{p}$ & Weighting assigning to timestep $p$\\
$CV_{l}$ & Conjectural variation associated with firm $l$\\
$E_{t}$ & Emissions factor level for technology $t$\\
 \hline
\end{tabular}%
\end{table}

In this section, we introduce the market we consider and describe the data inputs for models we use in Sections \ref{sec:results_mcp} and \ref{sec:model}. The electricity market we consider consists of two types of players: price-making firms and price-taking firms. Price-making firms may exert market power by using generation decisions to influence the market price. Price-taking firms do not have such ability. 

Each firm chooses its forward market generation decision so as to maximise its profits. Each firm may also hold multiple generating units with the technologies considered being baseload, mid merit and
peakload. The firms are distinguished by their price-making ability and their initial generation portfolio they may hold. However, each firm may also invest in new generation capacity in any of the technologies.

\begin{table}[htp!]
 \centering
\footnotesize
\caption{Initial power generation portfolio by firm ($CAP_{f,t}$). }\label{tab:firm_data}
\begin{tabular}{p{4cm}p{2cm}p{2cm}p{2cm}p{2cm}}
\hline
Technology & firm 1 & firm 2 & firm 3 & firm 4 \\
\hline
& price-making & price-making & price-taking & price-taking\\
\hline
Existing baseload (MW) & 1947 & 1940 & - & -  \\ 
Existing mid merit (MW) & 512 & - & 404 & - \\ 
Existing peakload (MW) & 270 & - & - & 234  \\ 
New baseload (MW) & 0 & 0 & 0 & 0  \\ 
New mid merit (MW) & 0 & 0 & 0 & 0 \\ 
New peakload (MW) & 0 & 0 & 0 & 0  \\ 
\hline
\end{tabular}%
\end{table}

We consider a electricity market that consists of four generating firms; firms $l=1$ and $l=2$ are price-making firms and while firms $f=3$ and $f=4$ are price-taking firms. We consider $|T|=6$ generating technologies; existing baseload, existing mid-merit, existing peaking, new baseload, new mid-merit and new peaking.  Each of the four firms hold different initial generating capacities. Firms $l=2$, $f=3$ and $f=4$ are, initially, specialised baseload, mid-merit and peaking firms, respectively. In contrast, firm $l=1$ is a integrated firm initially holding capacity across each of the existing technologies. Because of their sizes, the integrated firm and the specialised baseload firm are modelled as the price-making firms while the specialised mid-merit and peaking firms are the price-taking firms. Initially, each firm only holds `existing' technologies but, through their respective optimisation problems, may invest in any of the `new' technologies. Given the stylised nature of the model and following \cite{devine2019effect} and \cite{lynch2019impacts} , we do not explicitly model renewable technologies. Wind is incorporated into the model via the (net) demand intercept (see Market Clearing Condition \eqref{eqn:MCC}). We assume wind is not owned by any generation firm and its sole function is to reduce net demand. This is because wind has a marginal cost of zero and furthermore can only be dispatched downwards, and so given an exogenously-determined level of wind capacity, wind generation itself is unlikely to ever be strategically withheld by a generation firm \citep{devine2019effect}.

The initial portfolios of each firm are displayed in Table \ref{tab:firm_data}. These capacities follow from \cite{lynch2019role} and \cite{bertsch2018analysing} and are broadly based on \cite{eirgrid2016}.

\begin{table}[htp!]
 \centering
\footnotesize
\caption{Summary of techno-economic input data of considered supply side technologies. }\label{tab:data_tech}
\begin{tabular}{p{3cm}p{4cm}p{3cm}p{3cm}}
\hline
Technology & Annuity of specific invest  & Marginal gen. costs & Spec. CO$_2$ emissions \\
&($IC^{\text{GEN}}_{t})$ &($A^{\text{GEN}}_{t}$)& ($E_{t}$)\\ 
&(\euro/MW y) &(\euro/MWh$_{el}$)&(t CO$_2$/MWh$_{el}$)\\
\hline
Existing baseload & -  & 48.87 & 1.17 \\ 
Existing mid merit & -  & 41.10 & 0.36 \\ 
Existing peakload & -  & 63.38 & 0.56 \\ 
New baseload & 110,769  & 31.58 & 0.78 \\ 
New mid merit & 67,268  & 34.00 & 0.30 \\ 
New peakload & 40,363  & 50.50 & 0.45 \\ 
\hline
\end{tabular}%
\end{table}

The different characteristics associated with the technologies are displayed in Table \ref{tab:data_tech}. Both the marginal generation and investment costs again follow from \cite{lynch2019role} and \cite{bertsch2018analysing}. The marginal investment costs represent annualised investment costs. We consider $|P|=5$ forward time periods.  Table \ref{tab:demand} displays the demand curve intercept values which correspond to average hourly values for each time period.  However, each time period $p$ is assigned a weight $W_{p}=\frac{8760}{5}$. Thus, the test case in this work represents one year. Following \cite{lynch2017investment}, the five time periods represent summer low demand, summer high demand, winter low demand, winter high demand and winter peak demand. The demand curve slope value chosen is $B=9.091$. This parameter choice follows from \cite{devine2019effect} and \cite{di2013carbon}.

\begin{table}[htbp]
  \centering
    \begin{tabular}{rccccc}
    \hline
    Time Period ($p$) & 1     & 2     & 3     & 4     & 5 \\
    \hline
   &  25175.993   & 26768.307   & 30429.701   & 34302.196   & 37465.783 \\
    \hline
    \end{tabular}%
  \caption{Demand curve intercept ($A_{p}$) values}
  \label{tab:demand}%
\end{table}%

\section{Modelling electricity market as a mixed complementarity problem}\label{sec:results_mcp} 
In this section, we motivate our EPEC modelling approach, by applying the above data to a Mixed Complementarity Problem (MCP). This analysis demonstrates the naivety of using a MCP to model an oligopoly with a competitive fringe where both price-making and price-taking firms have investment decisions. A MCP determines an equilibrium of multiple optimisation problems by finding a point that satisfies the KKT conditions of each optimisation simultaneously  as a system of non-linear equations \citep{gabriel2009solving}. MCPs have used to model many energy markets \citep{huppmann2013endogenous, egging2013benders}. However, in a MCP modelling framework, a price-making firm's optimisation problem does not contain the optimal reactions of price-taking firms as constraints. As this subsection shows, this omission leads to myopic and contradictory outcomes. This is in contrast to the EPEC model described in Section \ref{sec:model}. 
\ref{sec:app_pm} describes the MCP problem. The MCP consists of the market clearing condition \eqref{eqn:MCC},  the price-taking firms' KKT conditions and the KKT conditions for all price-making firms. The parameter $CV_{l}$ represents the Conjectural Variation (CV) associated with firm $l$. CVs have been implemented in many cases MCP models \citep{egging2008complementarity, haftendorn2010modeling, huppmann2014market, egging2016risks} as they allow price-making firms to somewhat account for the optimal reactions of competitors. Conjectural variations take a value in the range $[0,1]$. \cite{huppmann2013endogenous} proposes a methodology to determine CVs that can be used to overcome myopic behaviour in models of an oligopoly with a competitive fringe. We advance the work of \cite{huppmann2013endogenous} by incorporating investment decisions.

We solve the MCP eleven times. Each time with a different CV for the two price-making firms; both firms have the same CV in each case. When $CV_{l=1}=CV_{l=2}=0$, both price-making firms lose their price-making ability and thus the market outcome corresponds to perfect competition. The remaining cases correspond to an oligopoly with competitive fringe modeled through conjectural variations. 

\begin{figure}[htbp]
	\centering
		\includegraphics[width=1\textwidth]{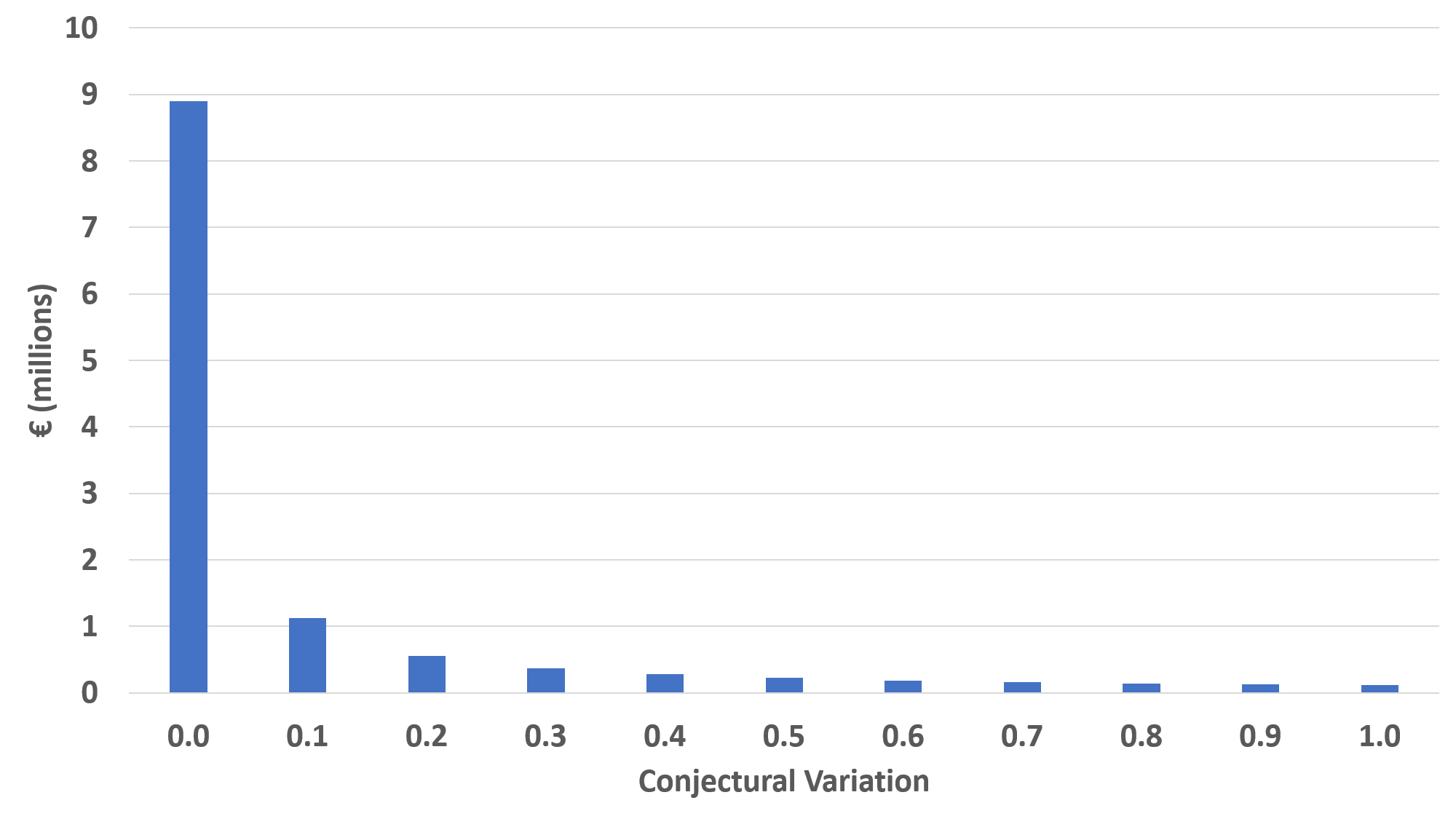}
		\caption{Price-making firm $l=1$'s profits using MCP setting framework.}\label{fig:myopic_profits}
\end{figure}

\begin{figure}[htbp]
	\centering
		\includegraphics[width=1\textwidth]{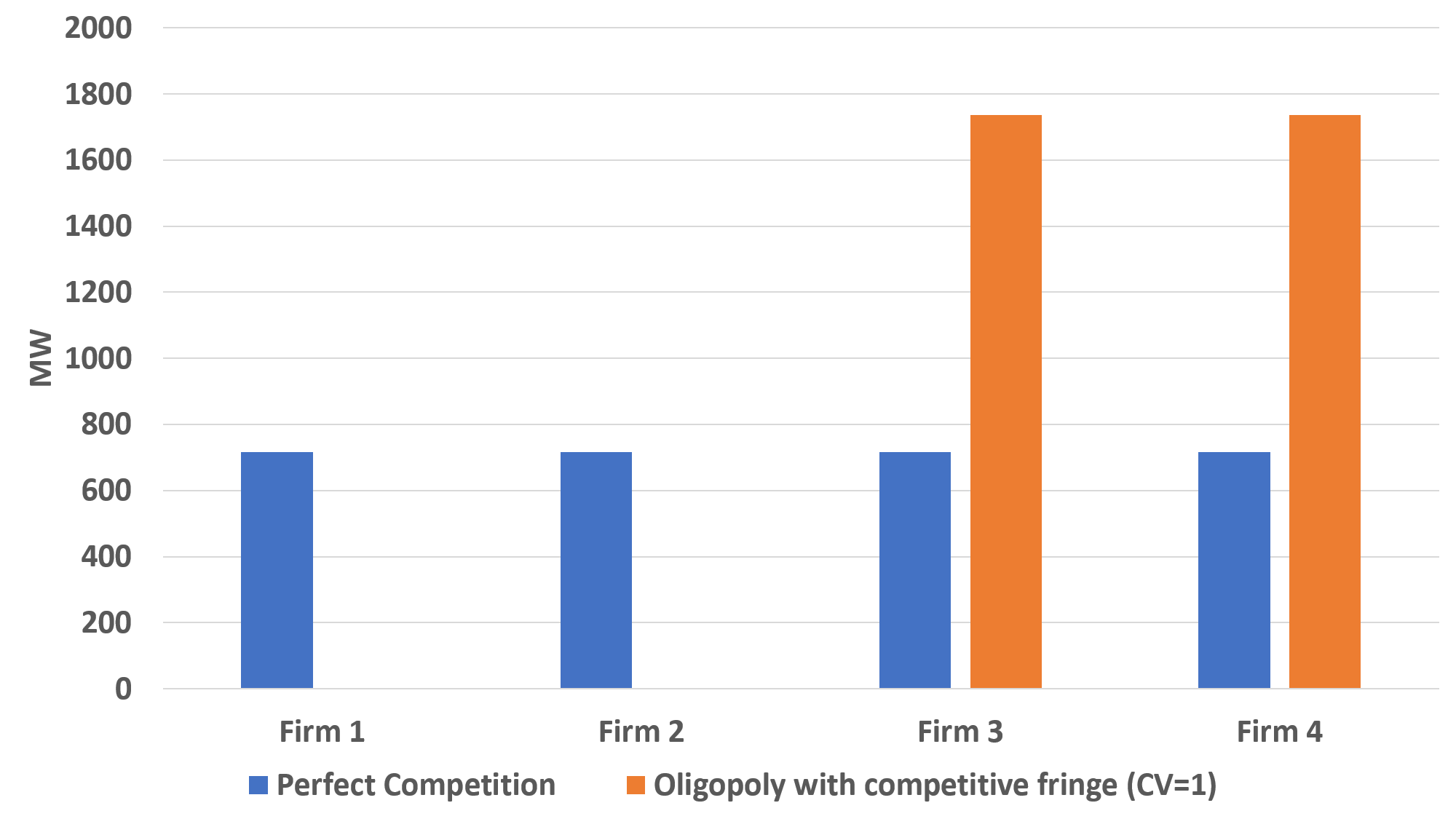}
		\caption{Price-making firms' investment into new mid-merit generation under MCP setting framework.}\label{fig:myopic_invest}
\end{figure}

Figure \ref{fig:myopic_profits} describes the profits of the price-making firm $l=1$ for the MCP cases. It shows that firm $l=1$ actually make less profits in the oligopoly with competitive fringe cases compared with perfect competition case. Clearly, if firms have price-making ability then they should be able to  to use that ability, at the very least, to make the same profits as they would have in a perfect competition setting. 

The result can be explained by Figure \ref{fig:myopic_invest} which shows the investment in new mid-merit generation for perfect competition case and $CV_{l=1}=CV_{l=2}=1$ case (similar results are observed for the $0<CV_{l}<1$ cases). In the perfect competition case, all firms invest 713MW into of new mid-merit generation.  However, in the oligopoly with competitive fringe case, the two price-making firms do not invest in any new technology while the two price-taking firms each increase their investment in new mid-merit generation to 1736MW. Note: in both cases, there are zero investments in new baseload or new peaking generation.  As equations \eqref{eqn:kkt_firm_gen_MCP} and \eqref{eqn:delta_gamma} show, price-making firms 1 and 2 assume  $\frac{\partial \gamma_{p}}{\partial  gen^{\text{PM}}_{l, t,p}}=-CV_{l}\times B$ in the oligopoly with competitive fringe case. These means that these two firms assume that if they decrease the amount of electricity they generate by one MW, then the equilibrium market price will increase by \euro $CV_{l}\times B$. In seeking to increase profits, price-making firms $l=1$ and $l=2$ decrease their generation in this way and hence do no invest in any new technology, as there is no point if they are not going to fully use that new generation. 

However, in the MCP with CV setting, price-making firms do not correctly account for the optimal reactions of the competitive fringe. Consequently, when price-making firms decrease their generation, the price-taking firms increase their generation and replace the price-making firms generation. Thus, the market price does not increase as much the price-making firms anticipate, if it increases at all. The expanded opportunity for price-taking firms to generate enables them to invest further into new mid-merit generation, as evidenced by Figure \ref{fig:myopic_invest}. 

The assumption that $\frac{\partial \gamma_{p}}{\partial  gen^{\text{PM}}_{l, t,p}}=-CV_{l}\times B$ in the oligopoly with competitive fringe case is not valid in a MCP setting where investment decisions are also incorporated. However, this assumption is valid in a MCP setting if all firms are price-making firms and hence behave in the same manner. In other words, when one firms seeks to increase the market price by decreases its generation, so does the rest of the firms and no one firm replaces the decreased generation from any other firm. 

This section demonstrates how the MCP modelling approach is unsuited to modelling an oligopoly with competitive fringe when investment decisions are also accounted for. Moreover, conjectural variations cannot overcome this modelling issue. In the following section, we show how the EPEC modelling approach overcomes the short-sighted/myopic behaviour observed in this section.

\section{Equilibrium problem with equilibrium constraints}\label{sec:model}
In this section we describe the Equilibrium Problem with Equilibrium Constraints (EPEC) we introduce in this work. As before, it represents an electricity market with two types of players: price-making firms and price-taking firms. Price-making firms may exert market power by using generation decisions to influence the market price. Price-taking firms do not have such ability. 

Each firm chooses its forward market generation decision so as to maximise its profits. Each firm may also hold multiple generating units with the technologies considered being baseload, mid merit and
peakload. The firms are distinguished by their price-making ability and their initial generation portfolio they may hold. However, each firm may also invest in new generation capacity in any of the technologies. 

The optimisation problems of each price-taking firm are embedded into the optimisation problem of each price-making firm. Thus, each price-making firm's problem is a bi-level optimisation problem and can be described as a Mathematical Program with Equilibrium Constraints (MPEC); the equilibrium constraints are the optimality conditions of the price-taking firms. This problem formulation enables each price-making firm to influence the market price through their decision variables, account for the optimal reactions of price-taking firms and, consequently, maximise profits. 

The overall EPEC problem is to find an equilibrium among the MPEC problems of each price-making firm which represent Nash Equilibria amongst them. Each MPEC problem can be represented as a Mixed Integer Non-Linear Problem (MINLP) and, thus, finding Nash Equilibria is a challenging task. To do so, we employ the Gauss-Seidel algorithm \cite{}. Furthermore to obtain an initial starting solution for this algorithm we utilise the approach taken in \cite{leyffer2010solving} for solving EPEC problems (henceforth known as the Leyffer-Munson approach). Both the Guass-Seidel algorithm and the Leyffer-Munson approach are described in detail in this section. 

Throughout this section the following conventions are used: lower-case Roman letters indicate indices or primal variables, upper-case Roman letters represent parameters (i.e., data), while Greek letters indicate prices or dual variables. The variables in parentheses alongside each constraint in this section are the Lagrange multipliers associated with those constraints. Tables \ref{tab:sets} - \ref{tab:parameters} explain the indices, variables and parameters, respectively, associated with both the price-making and price-taking firms' optimisation problems.

\subsection{Price-taking firm $f$'s problem}\label{sec:price_takers}
Price-taking firm $f$ seeks to maximise its profits (revenue less costs) by choosing investments in new capacity ($inv^{\text{PT}}_{f,t}$)  and by choosing the amount of electricity to generate from each technology at each time period ($gen^{\text{PT}}_{f, t,p}$). We assume each time period represents a forward time-period. Thus $gen^{\text{PT}}_{f, t,p}$ represents forward market generation decisions. Firm $f$'s costs include the per unit investment cost ($IC^{\text{GEN}}_{t}$) and the marginal cost of generation ($C^{\text{GEN}}_{t}$) while its revenues comes from the forward market price $\gamma_{p}$. 

Price-taking firm $f$'s optimisation problem is as follows:

\begin{equation}\label{eqn:Fgen_obj}
\begin{split}
\max_{\substack{ gen^{\text{PT}}_{f,t,p},
inv^{\text{PT}}_{f,t}\\
 }}\>\>
&\sum_{t,p} W_{p}\times gen^{\text{PT}}_{f, t,p} \times \big (  \gamma_{p}   -  C^{\text{GEN}}_{t}\big )-\sum_{t} IC^{\text{GEN}}_{t}\times inv^{\text{PT}}_{f,t},
\end{split}
\end{equation}
subject to:
\begin{equation}\label{eqn:Fgen_energy_con}
 gen^{\text{PT}}_{f, t,p} \leq CAP^{\text{PT}}_{f,t}+inv^{\text{PT}}_{f,t}, \>\> \forall t,p,
\end{equation}
where the parameter $W_{p}$ is the weight associated with timestep $p$. Constraint \eqref{eqn:Fgen_energy_con} ensures, for each generating technology in each timestep, firm $f$ cannot generate more than its initial capacity ($CAP^{\text{PT}}_{f,t}$) plus any new investments.  The variable alongside constraint \eqref{eqn:Fgen_energy_con} ($\lambda^{\text{PT},1}_{f,t,p}$) is the Lagrange multiplier associated with this constraint. In addition, each of firm $f$'s generation and investment decisions are constrained to be non-negative.

As firm $f$ is a price-taker, it cannot influence the market price with its generation decision.  The variable $\gamma_{p}$ is exogenous to firm $f$'s problem but a variable of the overall EPEC problem. When determining firm $f$'s Karush-Kuhn-Tucker (KKT) conditions, it is assumed $\frac{\partial \gamma_{p}}{\partial  gen^{\text{PT}}_{f, t,p}}=0$. Moreover, firm $f$ cannot see and hence account for the optimal decisions of other price-taking firms in addition to those of price-making firms. As firm $f$'s problem is linear, solving its associated KKT conditions ensures its problem is optimised \citep{gabriel2012complementarity}.

\subsection{Market clearing conditions}\label{sec:MCC}
The forward market price for each time period is determined from the following market clearing condition:
\begin{equation}\label{eqn:MCC}
\gamma_{p}=A_{p} -B \times \big( \sum_{ll,tt} gen^{\text{PM}}_{ll,tt,p}+\sum_{ff,tt} gen^{\text{PT}}_{ff,tt,p} \big), \>\> \forall p,
\end{equation}
where $A_{p}$ represents the demand curve intercept for each time period while $B$ is the time independent demand curve slope. Condition \eqref{eqn:MCC}  represents a linear demand curve and allows the market price to increases as the total market generation decreases and vice-versa.

\subsection{Price-maker $l$'s MPEC}
Price-making firm $l$'s optimisation problem is similar to firm $f$'s problem in that it too seeks to maximise profits (revenues less cost) by choosing its investment ($inv^{\text{PM}}_{l,t}$) and forward market generation ($gen^{\text{PM}}_{l, t,p}$) decisions. As before, firm $l$'s revenues come the forward market price while its costs include marginal generation and investment costs. In contrast to Section \ref{sec:price_takers} however, price-making firm $l$ can use their generation decisions to influence the market price. Price-making firm $l$ can also account for the optimal reactions of the price-taking firms. Its objective function is 
\begin{equation}\label{eqn:Lgen_obj1}
\begin{split}
\max_{\substack{ 
gen^{\text{PM}}_{l,t,p}, inv^{\text{PM}}_{l,t}\\
gen^{\text{PT}}_{f,t,p}, inv^{\text{PT}}_{f,t}\\
\gamma_{p}, \lambda^{\text{PT}}_{f, t,p}
 }}\>\>
&\sum_{t,p} W_{p}\times gen^{\text{PM}}_{l, t,p} \times \big (  \gamma_{p}   -  C^{\text{GEN}}_{t}\big )-\sum_{t} IC^{\text{GEN}}_{t}\times inv^{\text{PM}}_{l,t}.
\end{split}
\end{equation}
As firm $l$ can influence the market price through its generation decisions, we re-write objective function \eqref{eqn:Lgen_obj1} using market clearing condition \eqref{eqn:MCC} as follows:  
\begin{equation}\label{eqn:Lgen_obj2}
\begin{split}
\max_{\substack{ 
gen^{\text{PM}}_{l,t,p}, inv^{\text{PM}}_{l,t}\\
gen^{\text{PT}}_{f,t,p}, inv^{\text{PT}}_{f,t}\\
\gamma_{p}, \lambda^{\text{PT}}_{f, t,p}
 }}\>\>
&\sum_{t,p} W_{p} \times \bigg  (  A_{p} -B \times \big( \sum_{ll,tt} gen^{\text{PM}}_{ll,tt,p}+\sum_{ff,tt} gen^{\text{PT}}_{ff,tt,p} \big)      -  C^{\text{GEN}}_{t}  \bigg ) \times gen^{\text{PM}}_{l, t,p} -\sum_{t} IC^{\text{GEN}}_{t}\times inv^{\text{PM}}_{l,t}.
\end{split}
\end{equation}
The constraints of price-making firm $l$'s problem are
\begin{eqnarray}
 gen^{\text{PM}}_{l, t,p} &\leq& CAP^{\text{PM}}_{l,t}+inv^{\text{PM}}_{l,t}, \>\> \forall t,p, \label{eqn:Lgen_energy_con}
\end{eqnarray}
where each of firm $l$'s generation and investment decisions are also constrained to be non-negative. As with the price-taking firms, constraint \eqref{eqn:Lgen_energy_con} ensures, for each technology at each time period, firm $l$ cannot generate more electricity than its initial initial capacity plus any new investments. In addition to constraint \eqref{eqn:Lgen_energy_con}, firm $l$'s constraints also include the KKT conditions of the each price-taking firm:
\begin{eqnarray}
0\leq gen^{\text{PT}}_{f, t,p} &\perp& -W_{p} \times \big( \gamma_{p}  
-C^{\text{GEN}}_{t} \big)+\lambda^{\text{PT}}_{f, t,p} \geq 0, \>\> \forall f,t,p,\label{eqn:kkt_firm_gen}\\
0 \leq inv^{\text{PT}}_{f, t} &\perp& IC^{\text{GEN}}_{t}-\sum_{p}\lambda^{\text{PT}}_{f, t,p} \geq 0, \>\> \forall f,t,\label{eqn:kkt_firm_inv}\\
0 \leq \lambda^{\text{PT},1}_{f,t,p} &\perp& -gen^{\text{PT}}_{f, t,p}+CAP^{\text{PT}}_{f,t}+inv^{\text{PT}}_{f,t} \geq 0, \>\> \forall f,t,p.\label{eqn:kkt_firm_lambda}
\end{eqnarray}
Using market clearing condition \eqref{eqn:MCC} leads to condition \eqref{eqn:kkt_firm_gen} being re-written as follows:
\begin{equation}\label{eqn:kkt_firm_gen_comb}
    0\leq gen^{\text{PT}}_{f, t,p} \perp  -W_{p} \times \bigg(A_{p} -B \times \big( \sum_{ll,tt} gen^{\text{PM}}_{ll,tt,p}+\sum_{ff,tt} gen^{\text{PT}}_{ff,tt,p} \big) 
-C^{\text{GEN}}_{t} \bigg)+\lambda^{\text{PT}}_{f, t,p} \geq 0, \>\> \forall f,t,p.
\end{equation}
Constraints \eqref{eqn:kkt_firm_inv} - \eqref{eqn:kkt_firm_gen_comb} represent the optimal reactions, of each price-taking firm. As firm $f$'s problem (equations \eqref{eqn:Fgen_obj} and \eqref{eqn:Fgen_energy_con}) is a linear optimisation problem, these KKT conditions are both necessary and sufficient for optimality for the price-taking firms \citep{gabriel2012complementarity}.

Incorporating these conditions as constraints makes firm $l$'s problem a bi-level optimisation problem and ensures firm $l$ correctly anticipates how each price-taking firm will react to its decisions. Thus, this allows firm $l$ to adjust its decisions accordingly when seeking maximise its profits. 

Firm $l$'s optimisation problem is affected by the generation decisions of all other price-making firms; see objective function \eqref{eqn:Lgen_obj2} and constraint \eqref{eqn:kkt_firm_gen_comb}. However, when solving firm $l$'s problem we assume the decisions of all other price-making firms are fixed and exogenous to firm $l$'s problem. Sections \ref{sec:overall_EPEC} -- \ref{sec:overll_algorithm} describes how the optimisation problems all price-making firms are solved such that solutions represent Nash equilibria. 

As the KKT conditions \eqref{eqn:kkt_firm_inv} - \eqref{eqn:kkt_firm_gen_comb} represent the equilibrium constraints (optimal reactions) of the price-taking firms, firm $l$'s problem is a Mathematical Program with Equilibrium Constraints. We denote this problem as MPEC$_{l}$, which is a non-linear mathematical program because of the bi-linear terms in objective function \eqref{eqn:Lgen_obj2} involving firms $l$'s generation decisions and the generation decisions of all price-taking firms and because of the complementarity conditions incorporated into constraints \eqref{eqn:kkt_firm_inv} - \eqref{eqn:kkt_firm_gen_comb}. However, following the approach presented in \cite{fortuny1981representation}, we can remove the latter source of non-linearity using disjunctive constraints and big $M$ notation. Consequently, this leads to constraints \eqref{eqn:kkt_firm_inv} - \eqref{eqn:kkt_firm_gen_comb} being re-written as follows:
\begin{eqnarray}
0 &\leq & gen^{\text{PT}}_{f, t,p} \leq M \times r^{\text{1}}_{f, t,p},\>\> \forall f,t,p,\label{eqn:DC_firm_gen}\\
0 &\leq & -W_{p} \times \bigg(A_{p} -B \times \big( \sum_{ll,tt} gen^{\text{PM}}_{ll,tt,p}+\sum_{ff,tt} gen^{\text{PT}}_{ff,tt,p} \big)  
-C^{\text{GEN}}_{t} \bigg)+\lambda^{\text{PT}}_{f, t,p} \leq M \times(1-r^{\text{1}}_{f, t,p}),\>\> \forall f,t,p,\label{eqn:DC_kkt_firm_gen}\\
0 &\leq & inv^{\text{PT}}_{f, t} \leq M \times r^{\text{2}}_{f, t}, \>\> \forall f,t,\label{eqn:DC_firm_inv}\\
0 &\leq & IC^{\text{GEN}}_{t}-\sum_{p}\lambda^{\text{PT}}_{f, t,p} \leq M \times(1-r^{\text{2}}_{f, t} ), \>\> \forall f,t,\label{eqn:DC_kkt_firm_inv}\\
0 &\leq & \lambda^{\text{PT},1}_{f, t,p}  \leq M \times r^{\text{3}}_{f, t, p}, \>\> \forall f,t,p.\label{eqn:DC_firm_lambda}\\
0 &\leq & -gen^{\text{PT}}_{f, t,p}+CAP^{\text{PT}}_{f,t}+inv^{\text{PT}}_{f,t} \leq M \times(1-r^{\text{3}}_{f, t,p} ), \>\> \forall f,t,p.\label{eqn:DC_kkt_firm_lambda}
\end{eqnarray}
where $r^{\text{1}}_{f, t,p}$, $r^{\text{2}}_{f, t}$ and $r^{\text{3}}_{f, t,p}$ all represent binary 0-1 variables.

When solving the overall EPEC problem using the Gauss-Seidel algorithm (see Section \ref{sec:overall_EPEC}), MPEC$_{l}$ is characterized by objective function \eqref{eqn:Lgen_obj2}, subject to constraint \eqref{eqn:Lgen_energy_con} and constraints \eqref{eqn:DC_firm_gen} - \eqref{eqn:DC_kkt_firm_lambda}. Consequently price-making firm $l$'s optimisation problem is a Mixed Integer Non-Linear Problem. In Section \ref{sec:results_epec}, we use the DICOPT solver in GAMS to solve it.

\subsection{Overall EPEC}\label{sec:overall_EPEC}
The overall EPEC can be expressed as the problem of finding Nash equilibria among the price-makers $l$:
\begin{equation*}
\begin{split}
\text{Find:  } & \bigg\{inv^{\text{PM}}_{l=1,t},...,inv^{\text{PM}}_{l=L,t},\\
 &gen^{\text{PM}}_{l=1, t,p},...,gen^{\text{PM}}_{l=L, t,p},\\
 &inv^{\text{PT}}_{f=1,t},...,inv^{\text{PT}}_{f=F,t},\\
 &gen^{\text{PT}}_{f=1, t,p},...,gen^{\text{PT}}_{f=F, t,p},\\
 &\lambda^{\text{PT}}_{f=1,t,p},...,\lambda^{\text{PT}}_{f=F,t,p},\\
  & \gamma_{p} \bigg\} \text{ that solve:}\\
  & \text{MPEC}_{l} \text{ for each } l=1,...,L.
\end{split}
\end{equation*}

To find the such equilbiria, we implement the following Gauss-Seidel \citep{gabriel2012complementarity} algorithm. The algorithm iteratively solves each price-making firm's MPEC problem by fixing every other price-making firms' decisions, until it converges to a point where neither leader has an optimal deviation.

\begin{algorithm}[H]
\SetAlgoLined
 \While{$\sum_{l} | x_{l,k}-x_{l,k-1} | > TOL $ and $k < K$ }{
    \For{$l=1,..,L$}{
   Assume price maker$_{-l}$'s decision variables are fixed\;
   Solve MPEC$_{l}$\;
   }
 }
 \caption{Gauss-Seidel algorithm}\label{al:GS}
\end{algorithm}
where $TOL$ and $K$ represent a pre-defined convergence tolerance and a maximum number of allowable iterations, respectively. The vector $x_{l,t}$ represents the vector of all MPEC$_{l}$'s primal variables at iteration $k$.

\subsection{Obtaining an initial solution}\label{sec:LM}
To improve computational efficiency, we utilise the approach to obtaining a strong stationary point for EPECs, as described in \cite{leyffer2010solving}. We then use the statioanry point obtained as a starting point to the Gauss-Seidel algorithm. In this subsection, we describe how the Leyffer-Munson method as applied to the EPEC presented in this work. 


Firstly, we re-write price-making firm $l$'s problem, as defined by equations \eqref{eqn:Lgen_obj2} - \eqref{eqn:kkt_firm_gen_comb}, using slack variables. We do this by converting firm $l$'s inequality constraints into equality constraints as follows: 

\begin{equation}\label{eqn:Lgen_obj_LM}
\begin{split}
\max\>\>
&\sum_{t,p} W_{p} \times \bigg (  A_{p} -B \times \big( \sum_{ll,tt} gen^{\text{PM}}_{ll,tt,p}+\sum_{ff,tt} gen^{\text{PT}}_{ff,tt,p} \big)      -  C^{\text{GEN}}_{t}  \bigg ) \times gen^{\text{PM}}_{l, t,p} -\sum_{t} IC^{\text{GEN}}_{t}\times inv^{\text{PM}}_{l,t},
\end{split}
\end{equation}
subject to:
\begin{eqnarray}
 CAP^{\text{PM}}_{l,t}+inv^{\text{PM}}_{l,t}- gen^{\text{PM}}_{l, t,p}-s^{\text{CON\_LR}}_{l, t,p}&=&0,   \>\> \forall t,p, \>\> (\lambda^{PM}_{l, t,p}), \label{eqn:LM_lambda}\\
 gen^{\text{PM}}_{l, t,p} &\geq& 0, \>\> \forall t,p,\>\> (\chi^{\text{GEN}}_{l,t,p}), \label{eqn:chi_gen}\\
 inv^{\text{PM}}_{l, t} &\geq& 0, \>\> \forall t,\>\> (\chi^{\text{INV}}_{l,t}), \label{eqn:chi_inve}\\
s^{\text{CON\_LR}}_{l, t,p}&\geq& 0, \>\> \forall t,p,\>\> (\mu^{\text{CON\_LR}}_{l, t,p}), \label{eqn:mu_kkt_Lgen}\\
    -W_{p} \times \bigg (A_{p} -B \times \big( \sum_{ll,tt} gen^{\text{PM}}_{ll,tt,p}+\sum_{ff,tt} gen^{\text{PT}}_{ff,tt,p} \big)-C^{\text{GEN}}_{t} \bigg )+\lambda^{\text{PT}}_{f, t,p}-s^{\text{KKT\_GEN}}_{f, t,p} &=& 0, \>\> \forall f,t,p,\>\> (\alpha^{\text{KKT\_GEN}}_{l, f, t,p}), \label{eqn:s_kkt_gen}\\
gen^{\text{PT}}_{f, t,p} &\geq& 0, \>\> \forall f,t,p,\>\> (\mu^{\text{KKT\_GEN}}_{l, f, t,p}), \label{eqn:mu_kkt_gen}\\
 s^{\text{KKT\_GEN}}_{f, t,p} &\geq& 0, \>\> \forall f,t,p,\>\> (\mu^{\text{s\_KKT\_GEN}}_{l, f, t,p}), \label{eqn:mu_s_kkt_gen}\\
 gen^{\text{PT}}_{f, t,p} \times s^{\text{KKT\_GEN}}_{f, t,p} &=& 0, \>\> \forall f,t,p,\>\> (\mu^{\text{GEN\_s\_KKT\_GEN}}_{l, f, t,p}), \label{eqn:mu_gen_s_kkt_gen}\\
IC^{\text{GEN}}_{t}-\sum_{p}\lambda^{\text{PT}}_{f, t,p}-s^{\text{KKT\_INV}}_{f, t} &=& 0, \>\> \forall f,t,\>\> (\alpha^{\text{KKT\_INV}}_{l, f, t}), \label{eqn:s_kkt_inv}\\
inv^{\text{PT}}_{f, t} &\geq& 0, \>\> \forall f,t,\>\> (\mu^{\text{KKT\_INV}}_{l, f, t}), \label{eqn:mu_kkt_inv}\\
s^{\text{KKT\_INV}}_{f, t} &\geq& 0, \>\> \forall f,t,\>\> (\mu^{\text{KKT\_s\_INV}}_{l, f, t}), \label{eqn:mu_s_kkt_inv}\\
inv^{\text{PT}}_{f, t} \times s^{\text{KKT\_INV}}_{f, t} &=& 0, \>\> \forall f,t,\>\> (\mu^{\text{INV\_s\_KKT\_INV}}_{l, f, t}), \label{eqn:mu_inv_s_kkt_inv}\\
-gen^{\text{PT}}_{f, t,p}+CAP^{\text{PT}}_{f,t}+inv^{\text{PT}}_{f,t} - s^{\text{CON\_FR}}_{f, t,p} &=& 0\>\> \forall f,t,p, \>\> (\alpha^{\text{CON}}_{f, t, p}), \label{eqn:alpha_con}\\
\lambda^{\text{PT}}_{f, t,p} &\geq& 0, \>\> \forall f,t,p,\>\> (\mu^{\text{CON}}_{l, f, t,p}), \label{eqn:mu_con}\\
s^{\text{CON\_FR}}_{f, t,p} &\geq& 0, \>\> \forall f,t,p,\>\> (\mu^{\text{s\_CON}}_{l, f, t,p}), \label{eqn:mu_s_con}\\
 \lambda^{\text{PT}}_{f, t,p} \times s^{\text{CON\_FR}}_{f, t,p} &=& 0, \>\> \forall f,t,p,\>\> (\mu^{\text{CON\_s\_CON}}_{l, f, t,p}). \label{eqn:mu_CON_s_CON}
\end{eqnarray}
The variables in brackets alongside each of these constraints are the Lagrange multipliers associated with those constraints. Note: each of the multipliers has a subscript $l$ associated with it showing  how there are unique multipliers for each of the price-making firms problems.

Secondly, we find the stationary KKT conditions of the optimisation problem \eqref{eqn:Lgen_obj_LM} - \eqref{eqn:mu_CON_s_CON}. Let $\mathcal{L}_{l}$ be the Lagrangian associated with that problem.
\begin{equation}\label{eqn:LM_genL}
\begin{split}
&\\ \frac{\partial \mathcal{L}_{l}}{\partial gen^{\text{PM}}_{l,t,p}}:- W_{p} \times \bigg (A_{p} -B \times \big( \sum_{ll,t} gen^{\text{PM}}_{ll,tt,p}+\sum_{ff,tt} gen^{\text{PT}}_{ff,tt,p} \big)-B \times gen^{\text{PM}}_{l,t,p}
&\\-C^{\text{GEN}}_{t} \bigg)+\lambda^{\text{PT}}_{f, t,p}+\sum_{ff,tt}W_{p} \times B \times \alpha^{\text{KKT\_GEN}}_{l, ff, tt,p} - \chi^{\text{GEN}}_{l,t,p}=0, \>\> \forall t,p,
\end{split}
\end{equation}
\begin{equation}\label{eqn:LM_invL}
\frac{\partial \mathcal{L}_{l}}{\partial inv^{\text{PM}}_{l,t}}: IC^{\text{GEN}}_{t}-\sum_{p}\lambda^{\text{PM}}_{l, t,p}- \chi^{\text{INV}}_{l,t} =0, \>\> \forall t,\\
\end{equation}
\begin{equation}\label{eqn:LM_genF}
\begin{split}
&\\\frac{\partial \mathcal{L}_{l}}{\partial gen^{\text{PT}}_{f,t,p}}: 
\sum_{tt}W_{p} \times B \times gen^{\text{PM}}_{l,tt,p} + \sum_{ff,tt}W_{p} \times B \times \alpha^{\text{KKT\_GEN}}_{l, ff, tt,p} - \mu^{\text{KKT\_GEN}}_{l, f, t,p} 
&\\+ s^{\text{KKT\_GEN}}_{f, t,p} \times \mu^{\text{GEN\_s\_KKT\_GEN}}_{l, f, t,p} -\alpha^{\text{CON}}_{f, t, p} =0, \>\> \forall f,t,p,\\
\end{split}
\end{equation}

\begin{equation}\label{eqn:LM_invF}
\frac{\partial \mathcal{L}_{l}}{\partial inv^{\text{PT}}_{f,t}}:-\mu^{\text{KKT\_INV}}_{l, f, t}+s^{\text{KKT\_INV}}_{f, t} \times \mu^{\text{INV\_s\_KKT\_INV}}_{l, f, t}-\sum_{p}\alpha^{\text{CON}}_{f, t, p}=0, \>\> \forall f,t,\\
\end{equation}

\begin{equation}\label{eqn:LM_lambdaF}
    \frac{\partial \mathcal{L}_{l}}{\partial \lambda^{\text{PT}}_{f, t,p}}:  -\alpha^{\text{KKT\_GEN}}_{l, f, t,p}+\alpha^{\text{KKT\_INV}}_{l, f, t}-\mu^{\text{CON}}_{l, f, t,p}+s^{\text{CON\_FR}}_{f, t,p} \times \mu^{\text{CON\_s\_CON}}_{l, f, t,p}=0,\>\> \forall f,t,p,\\
\end{equation}

\begin{equation}\label{eqn:LM_s_conL}
    \frac{\partial \mathcal{L}_{l}}{\partial s^{\text{CON\_LR}}_{l, t,p}}: \lambda^{PM}_{l, t,p}-\mu^{\text{CON\_LR}}_{l, t,p}=0,\>\> \forall t,p,\\
\end{equation}

\begin{equation}\label{eqn:LM_s_kkt_gen}
  \frac{\partial \mathcal{L}_{l}}{\partial s^{\text{KKT\_GEN}}_{f, t,p}}: -\alpha^{\text{KKT\_GEN}}_{l, f, t,p}-\mu^{\text{s\_KKT\_GEN}}_{l, f, t,p}+ gen^{\text{PT}}_{f, t,p} \times \mu^{\text{GEN\_s\_KKT\_GEN}}_{l, f, t,p}=0,\>\> \forall f,t,p,\\
\end{equation}

\begin{equation}\label{eqn:LM_s_kkt_inv}
  \frac{\partial \mathcal{L}_{l}}{\partial s^{\text{KKT\_INV}}_{f, t}}: -\alpha^{\text{KKT\_INV}}_{l, f, t}-\mu^{\text{KKT\_s\_INV}}_{l, f, t}+inv^{\text{PT}}_{f, t} \times \mu^{\text{INV\_s\_KKT\_INV}}_{l, f, t}=0,\>\> \forall f,t,\\
\end{equation}

\begin{equation}\label{eqn:LM_s_conF}
  \frac{\partial \mathcal{L}_{l}}{\partial s^{\text{CON\_FR}}_{f, t,p}}: -\alpha^{\text{CON}}_{f, t, p}-\mu^{\text{s\_CON}}_{l, f, t,p}+\lambda^{\text{PT}}_{f, t,p} \times \mu^{\text{CON\_s\_CON}}_{l, f, t,p}=0,\>\> \forall f,t,p.\\
\end{equation}
In addition, each of the Lagrange multipliers associated with inequality constraints in \eqref{eqn:Lgen_obj_LM} - \eqref{eqn:mu_CON_s_CON} are constrained to be non-negative.

Following this, we find the complementary KKT conditions of the optimisation problem \eqref{eqn:Lgen_obj_LM} - \eqref{eqn:mu_CON_s_CON} as follows:

\begin{eqnarray}
 gen^{\text{PM}}_{l,t,p} \times \chi^{\text{GEN}}_{l,t,p} &=&0, \>\> \forall t,p,\label{eqn:LM_C_genL}\\
 inv^{\text{PM}}_{l,t} \times \chi^{\text{INV}}_{l,t} &=&0, \>\> \forall t,\label{eqn:LM_C_invL}\\
 s^{\text{CON\_LR}}_{l, t,p} \times \mu^{\text{CON\_LR}}_{l, t,p}&=&0, \>\> \forall t,p,\label{eqn:LM_C_s_conL}\\
gen^{\text{PT}}_{f, t,p} \times \mu^{\text{KKT\_GEN}}_{l, f, t,p}&=&0, \>\> \forall f,t,p,\label{eqn:LM_C_genF}\\
 s^{\text{KKT\_GEN}}_{f, t,p} \times \mu^{\text{s\_KKT\_GEN}}_{l, f, t,p}&=&0, \>\> \forall f,t,p,\label{eqn:LM_C_s_KKT_gen}\\
inv^{\text{PT}}_{f, t} \times \mu^{\text{KKT\_INV}}_{l, f, t}&=&0, \>\> \forall f,t,\label{eqn:LM_C_invF}\\
s^{\text{KKT\_INV}}_{f, t} \times \mu^{\text{KKT\_s\_INV}}_{l, f, t}&=&0, \>\> \forall f,t,\label{eqn:LM_C_s_KKT_inv}\\
\lambda^{\text{PT}}_{f, t,p} \times \mu^{\text{CON}}_{l, f, t,p}&=&0, \>\> \forall f,t,p,\label{eqn:LM_C_lambdaF}\\
s^{\text{CON\_FR}}_{f, t,p} \times \mu^{\text{s\_CON}}_{l, f, t,p}&=&0, \>\> \forall f,t,p.\label{eqn:LM_C_s_conF}
\end{eqnarray}

The Leyffer-Munson method, as applied to this work, is obtain a solution set that satisfies conditions \eqref{eqn:LM_lambda} - \eqref{eqn:LM_C_s_conF} of each price-making firm $l$ simultaneously. To do this, each KKT condition with bi-linear terms (equations \eqref{eqn:mu_gen_s_kkt_gen}, \eqref{eqn:mu_inv_s_kkt_inv}, \eqref{eqn:mu_CON_s_CON} and \eqref{eqn:LM_C_genL} - \eqref{eqn:LM_C_s_conF}) are removed as constraints and are summed together to create the following objective function:
\begin{equation}\label{eqn:LM_Full_Obj}
\begin{split}
\min &  \sum_{f,t,p} gen^{\text{PT}}_{f, t,p} \times s^{\text{KKT\_GEN}}_{f, t,p} +
 \sum_{f,t}inv^{\text{PT}}_{f, t} \times s^{\text{KKT\_INV}}_{f, t} +
 \sum_{f,t,p} \lambda^{\text{PT}}_{f, t,p} \times s^{\text{CON\_FR}}_{f, t,p} +
\sum_{l,t,p} gen^{\text{PM}}_{l,t,p} \times \chi^{\text{GEN}}_{l,t,p}+ 
\\& \sum_{l,t} inv^{\text{PM}}_{l,t} \times \chi^{\text{INV}}_{l,t} +
\sum_{l,t,p} s^{\text{CON\_LR}}_{l, t,p} \times \mu^{\text{CON\_LR}}_{l, t,p}+
\sum_{f,t,p} gen^{\text{PT}}_{f, t,p} \times \mu^{\text{KKT\_GEN}}_{l, f, t,p}+
\\& \sum_{f,t,p} s^{\text{KKT\_GEN}}_{f, t,p} \times \mu^{\text{s\_KKT\_GEN}}_{l, f, t,p}+
 \sum_{f,t} inv^{\text{PT}}_{f, t} \times \mu^{\text{KKT\_INV}}_{l, f, t}+
\sum_{f,t} s^{\text{KKT\_INV}}_{f, t} \times \mu^{\text{KKT\_s\_INV}}_{l, f, t}+
\\& \sum_{f,t,p} \lambda^{\text{PT}}_{f, t,p} \times \mu^{\text{CON}}_{l, f, t,p}+
\sum_{f,t,p} s^{\text{CON\_FR}}_{f, t,p} \times \mu^{\text{s\_CON}}_{l, f, t,p}.
\end{split}
\end{equation}
Thus, the Leyffer-Munson optimisation problem is to minimise equation \eqref{eqn:LM_Full_Obj} subject to constraints \eqref{eqn:LM_lambda} - \eqref{eqn:mu_s_kkt_gen}, \eqref{eqn:s_kkt_inv} - \eqref{eqn:mu_s_kkt_inv}, \eqref{eqn:alpha_con} - \eqref{eqn:mu_s_con}  and \eqref{eqn:LM_genL} - \eqref{eqn:LM_s_conF}. In addition, each of the Lagrange multipliers associated with inequality constraints in \eqref{eqn:Lgen_obj_LM} - \eqref{eqn:mu_CON_s_CON} are constrained to be non-negative. The Leyffer-Munson optimization problem is a Non-Linear Program (NLP) and, in Section \ref{sec:results_epec}, we use the CONOPT solver in GAMS to solve it.

\subsection{Overall algorithm}\label{sec:overll_algorithm}
Algorithm \ref{al:overall} describes the overall algorithm for finding Nash equilibria from the EPEC problem. For iteration $i$, we firstly provide a random solution set from the search space and use these as initial starting point solutions for the Leyffer-Munson approach. As the Leyffer-Munson approach is a non-linear optimisation problem, the CONOPT solver does not always find a local minimum. If the Leyffer-Munson method does not converge to a locally optimal solution, then iteration $i$ is deemed unsuccessful and the algorithm skips ahead to iteration $i+1$. If the Leyffer-Munson approach does converge however, the locally optimal solution is then used as starting point solution for the Gauss-Seidel algorithm. If the Gauss-Seidel does (not) converge to a Nash equilibrium solution, then iteration $i$ is (not) deemed successful.  This process is repeated for $I$ iterations.  
\begin{algorithm}[H]
\SetAlgoLined
 \For{$i=1,..,I$}{
 Provide random initial solutions\;
 Solve Leyffer-Munson optimisation problem\;
    \If{Solution from LM is locally optimal}{
  Solve EPEC using Gauss-Seidel algorithm using solutions from LM as starting point \;
  \If{Gauss-Seidel algorithm converges}{
  Save solution\;
  }
   }
 }
\caption{Overall algorithm for finding Nash Equilibria}\label{al:overall}
\end{algorithm}

\section{Results from EPEC model}\label{sec:results_epec}
In this section, we present the results when the data presented in Section \ref{sec:input_data} is applied to the EPEC model described in Section \ref{sec:model}. We focus on the firms' profits, firms' investment decisions, market prices, consumer costs and carbon emissions. To obtain these results we utilise the algorithm described in Section \ref{sec:overll_algorithm} for $I=2000$ iterations. For the first 1000 iterations firm $l=1$'s MPEC problem is solved before firm $l=2$'s MPEC problem. For the subsequent 1000 iterations the opposite applies and firm $l=2$'s MPEC problem is solved before firm $l=1$'s MPEC problem. 

The algorithm did not always find a Nash Equilibrium (NE) solution. In fact, in the results to follow, only 72 of the 2000 iterations successfully found a NE solution, henceforth know as successful iterations. Of these, 62 iterations occurred when firm $l=1$'s MPEC problem was solved before firm $l=2$'s MPEC problem while the remaining 10 successful iterations occurred when firm $l=2$'s MPEC problem was solved first. For unsuccessful iterations the algorithm failed to find a NE solution for one of two reasons:
\begin{enumerate}
    \item For the random initial solution provided, the Leyffer Munson was found to be locally infeasible by the CONOPT solver.
    \item For the Gauss-Seidel algorithm, the convergence tolerance remained greater than $TOL=10^{-3}$ after $K=100$ iterations.  
\end{enumerate}

For each of the 72 successful iterations, the Leyffer-Munson approach solved to a locally optimal solution which implied that it is not necessarily a feasible solution to the EPEC. For 43\% of these successful iterations, the objective function for the Leyffer-Munson approach (equation \eqref{eqn:LM_Full_Obj}) converged to zero. For the remaining 57\% of successful iterations, the objective function converged to a strictly positive objective function value.  When the Leyffer-Munson approach gives a non-zero objective function value, the solutions cannot guaranteed to be a feasible point for the overall EPEC. However, the results in this section show that, despite this, such solutions can still provide good starting point solutions to the Gauss-Seidel algorithm.

\begin{figure}[htbp]
	\centering
		\includegraphics[width=1\textwidth]{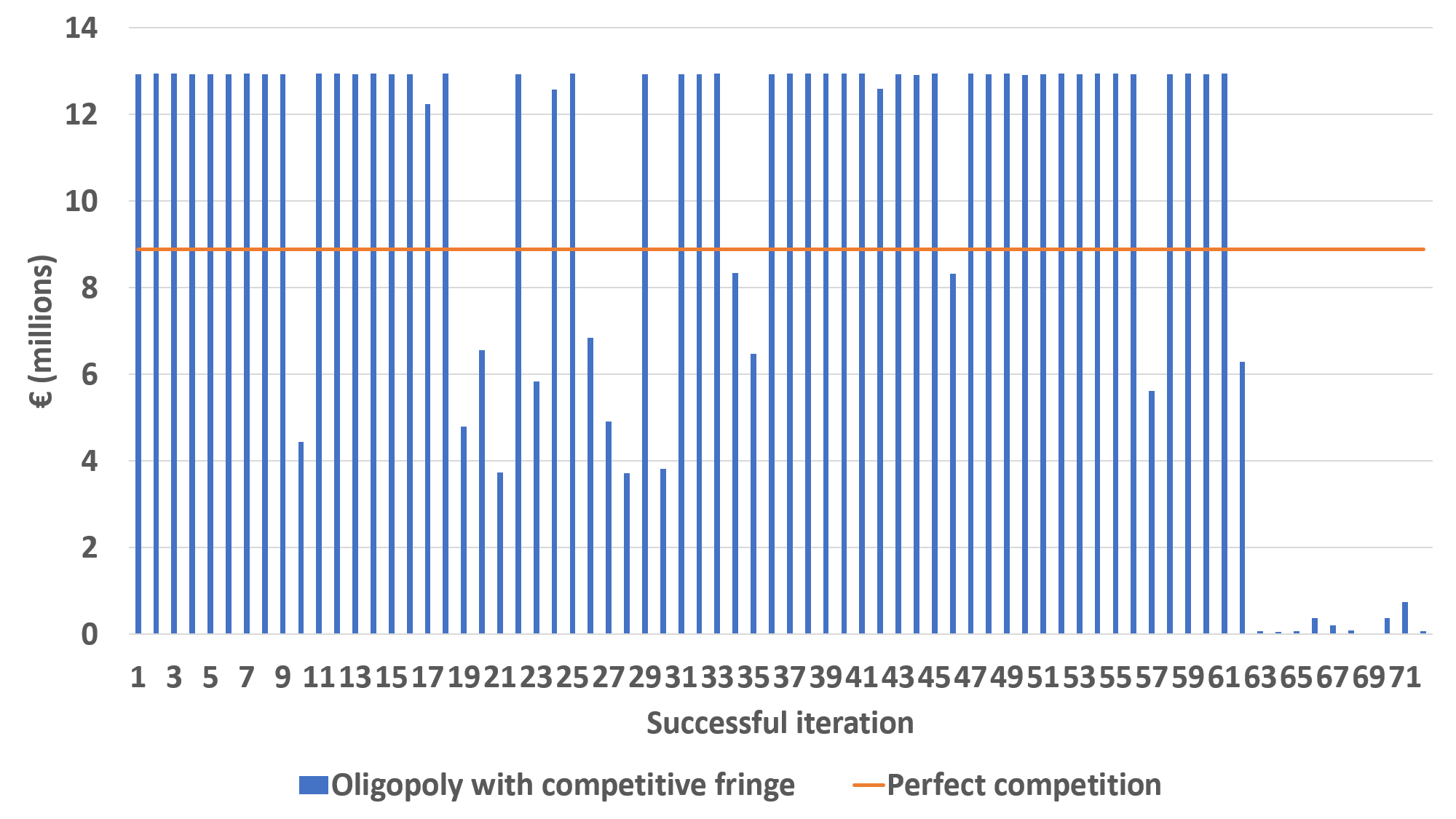}
		\caption{Profits for price-making firm $l=1$ for each successful iteration.}\label{fig:profits_l1}
\end{figure}

\begin{figure}[htbp]
	\centering
		\includegraphics[width=1\textwidth]{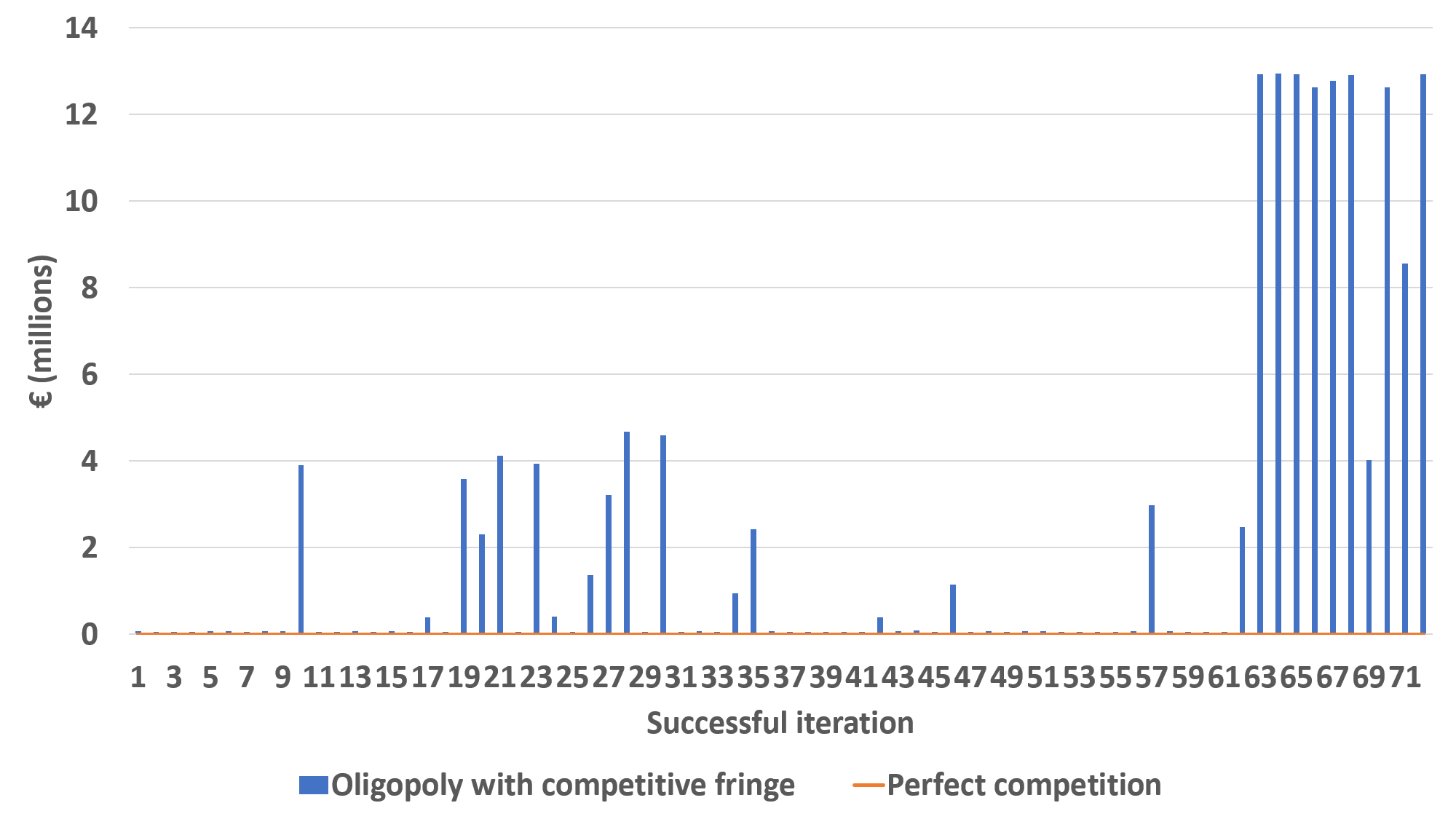}
			\caption{Profits for price-making firm $l=2$ for each successful iteration.}\label{fig:profits_l2}
\end{figure}

\begin{figure}[htbp]
	\centering
		\includegraphics[width=1\textwidth]{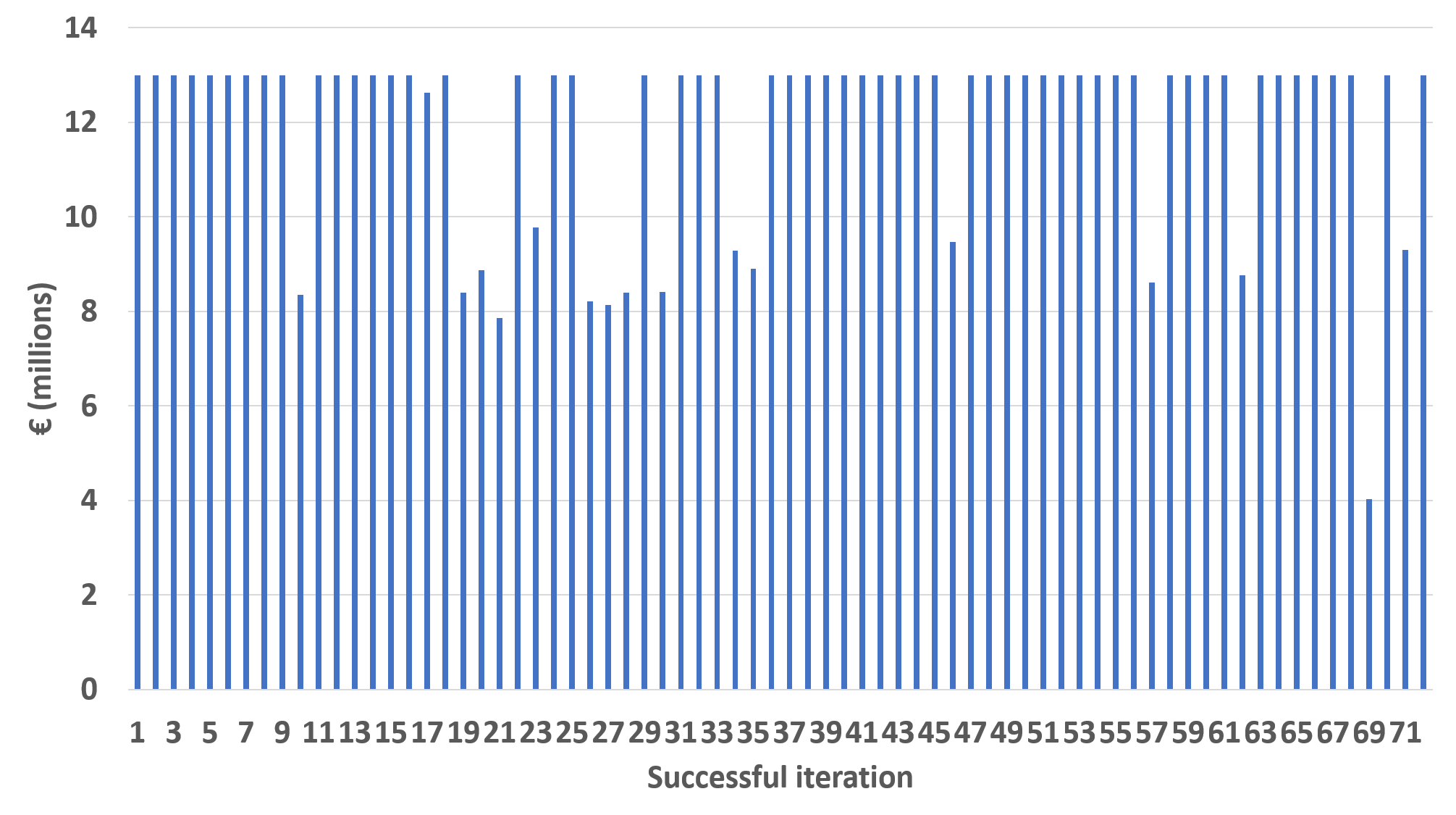}
			\caption{Combined profits for price-making firms for each successful iteration.}\label{fig:profits2}
\end{figure}

Figures \ref{fig:profits_l1} and \ref{fig:profits_l2} display price-making firms $l=1$ and $l=2$ profits, respectively, for each of the successful iterations. The horizontal lines in each figure represent the profits each firm would make from the perfect competition case in the Section \ref{sec:results_mcp}. Figures \ref{fig:profits_l1} and \ref{fig:profits_l2} both show that the algorithm found multiple NE solutions. Firm $l=1$'s profits varies from \euro 0 to \euro 12.98M while firm $l=2$'s profits ranged from \euro 50,000 to \euro 6.8M. For the majority of equilibria found, both firms made profits greater than they would in a perfect framework. This occurred in 66.7\% and 100\% of the successful iterations for firms $l=1$ and $l=2$ profits, respectively. Thus showing that there was no NE found where both firms' profits were below their perfect competition equivalent.  

Figure \ref{fig:profits2} display the combined profits of the two price-making firms. It shows that the combined profits varied across the equilibria suggesting that there was not a zero-sum game between the price-making firms on how profits were split between them. We describe below why there are multiple equilibria and why, in some equilibria, one of the price-making firms makes a profit less than it would in a perfectly competitive market. 

\begin{figure}[htbp]
	\centering
		\includegraphics[width=1\textwidth]{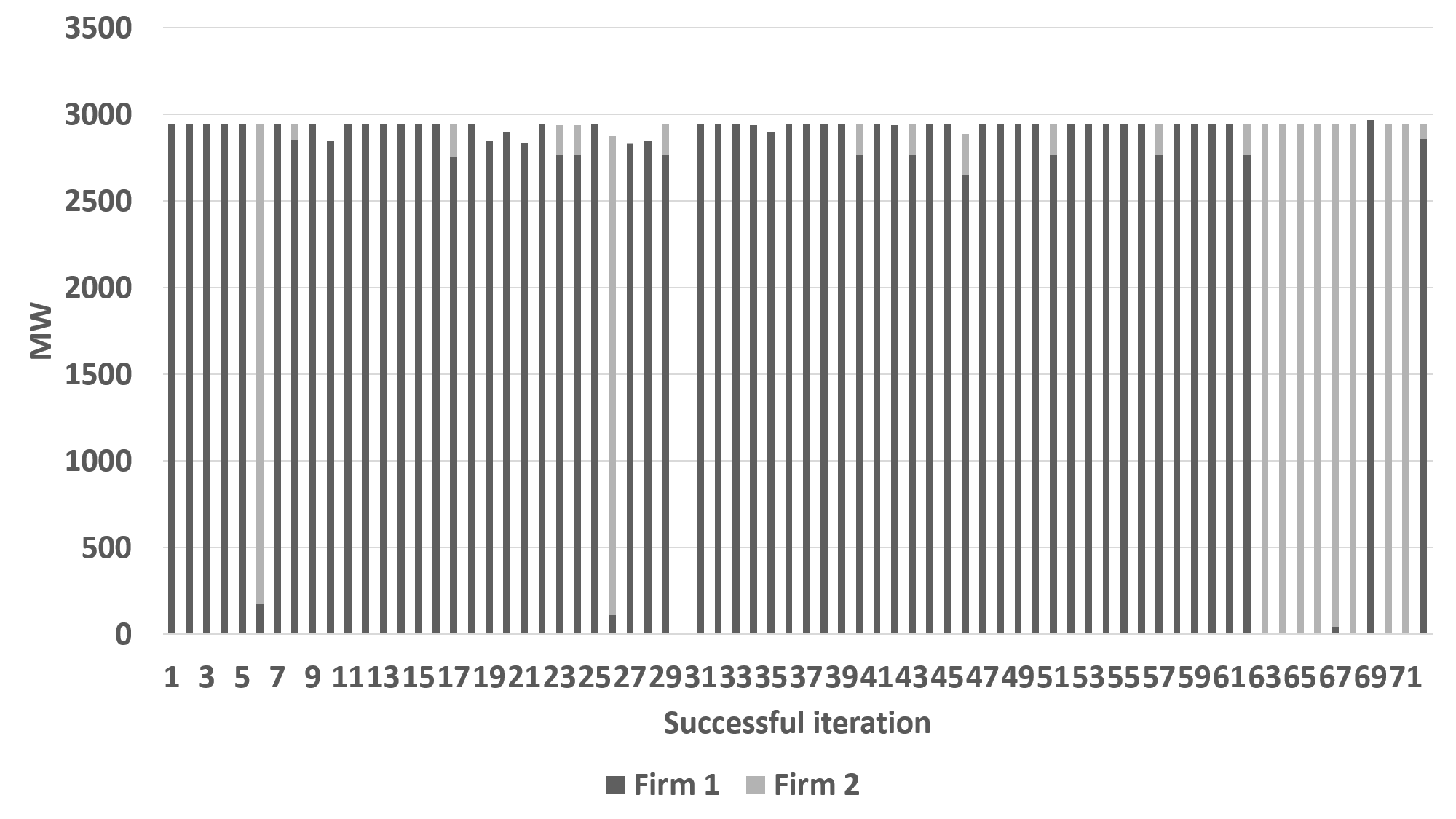}
			\caption{Combined investment in new mid-merit for price-making firms for each successful iteration.}\label{fig:invest}
\end{figure}

Figure \ref{fig:invest} displays the combined investments into new mid-merit generation for the two price-making firms for each successful iteration. In comparison to Section \ref{sec:results_mcp}, both of these firms did not invest in any baseload or peaking generation. However, in contrast to Section \ref{sec:results_mcp}, there was only one equilibrium point found where the price-taking firms invested in any new generation technology.

For the majority of equilibria found (80\%), the combined investments were 2941MW. But, for some equilibria, the combined investments were slightly higher with maximum combined investment reaching 2967MW at one equilibrium point while the lowest combined investment was 2831MW. The first 62 successful iterations in Figures \ref{fig:profits_l1} - \ref{fig:profits2} came when firm $l=1$'s MPEC problem was solved before firm $l=2$'s. At these equilibria, firm $l=1$ and $l=2$'s investments in new mid-merit generation averaged at 2807MW and 496MW, respectively. The final 10 successful iterations occurred when firm $l=2$'s MPEC problem was solved first. As a result, at these equilibria, firm $l=1$ average investments in new generation decreased to 1467MW while firm $l=2$'s increased to 2618MW. 

The results show that firm $l=1$ makes a profit less than it would of in a perfectly competitive market, in some of the equilibria found. This is because we model two price-making firms. When price-making firm $l$ commits to a large amount of forward generation, it can leave firm $\hat{l} \neq l$ with a reduced opportunity to generate and hence reduced profits. We explain this in further detail in Section \ref{sec:reasons}. Interestingly, in each of the successful iterations where firm $l=2$ MPEC is solved first, firm $l=1$'s profits are significantly below their perfect competition result while firm $l=2$ are significantly higher - this result highlights the importance of the order in which MPEC problems of an EPEC are searching, when searching for equilibria. 

\begin{figure}[htbp]
	\centering
		\includegraphics[width=1\textwidth]{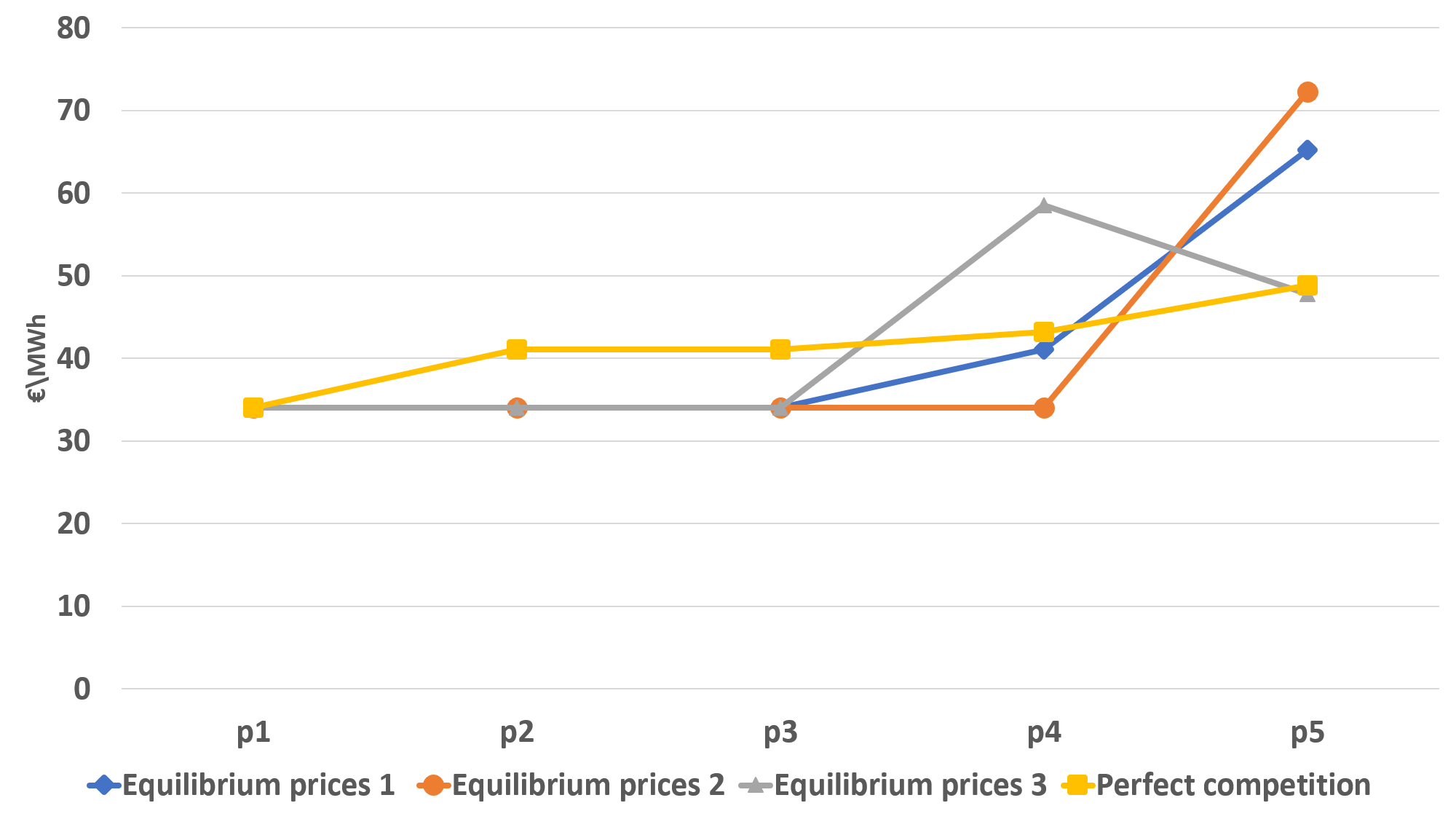}
			\caption{Equilibrium forward market prices ($\gamma_{p}$) for EPEC model versus perfect competition.}\label{fig:pices2}
\end{figure}

Despite Figures \ref{fig:profits_l1} - \ref{fig:invest} presenting the multiple equilibria, the forward prices rested at one of three price time series. Time series one and two were observed in 81.9\% and 16.7\% of the equilbira found while the third series was only observed at one of the equilibrium points found. Figure \ref{fig:pices2} displays these three price time series along with the prices from the perfect competition case of Section \ref{sec:results_mcp}. Interestingly, for time periods $p=2,3$, the forward market prices in the oligopoly with competitive fringe case are less than those from the perfect competition case. This is despite half the firms having price-making ability. However, the market prices in the oligopoly with competitive fringe case are higher at later time steps. Note: both the first and second equilibrium price time series were found when both firm $l=1$ and firm $l=2$ MPEC's were solved first while the only instance of the third equilibrium price time series occurred when firm $l=2$'s MPEC was solved first.

\begin{figure}[htbp]
	\centering
		\includegraphics[width=1\textwidth]{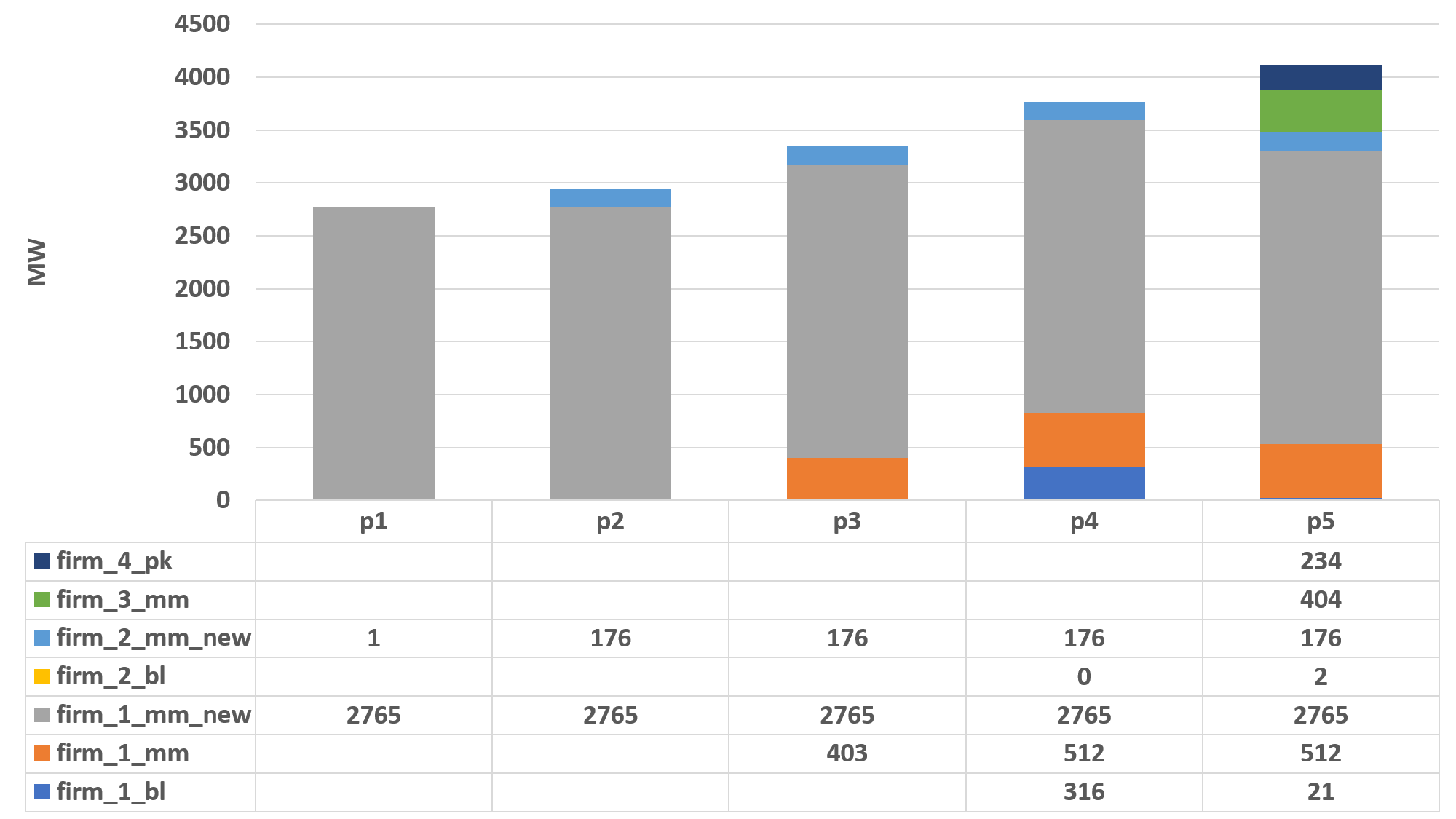}
		\caption{Generation mix for the first successful iteration.}\label{fig:gen_SI1}
\end{figure}

\begin{figure}[htbp]
	\centering
		\includegraphics[width=1\textwidth]{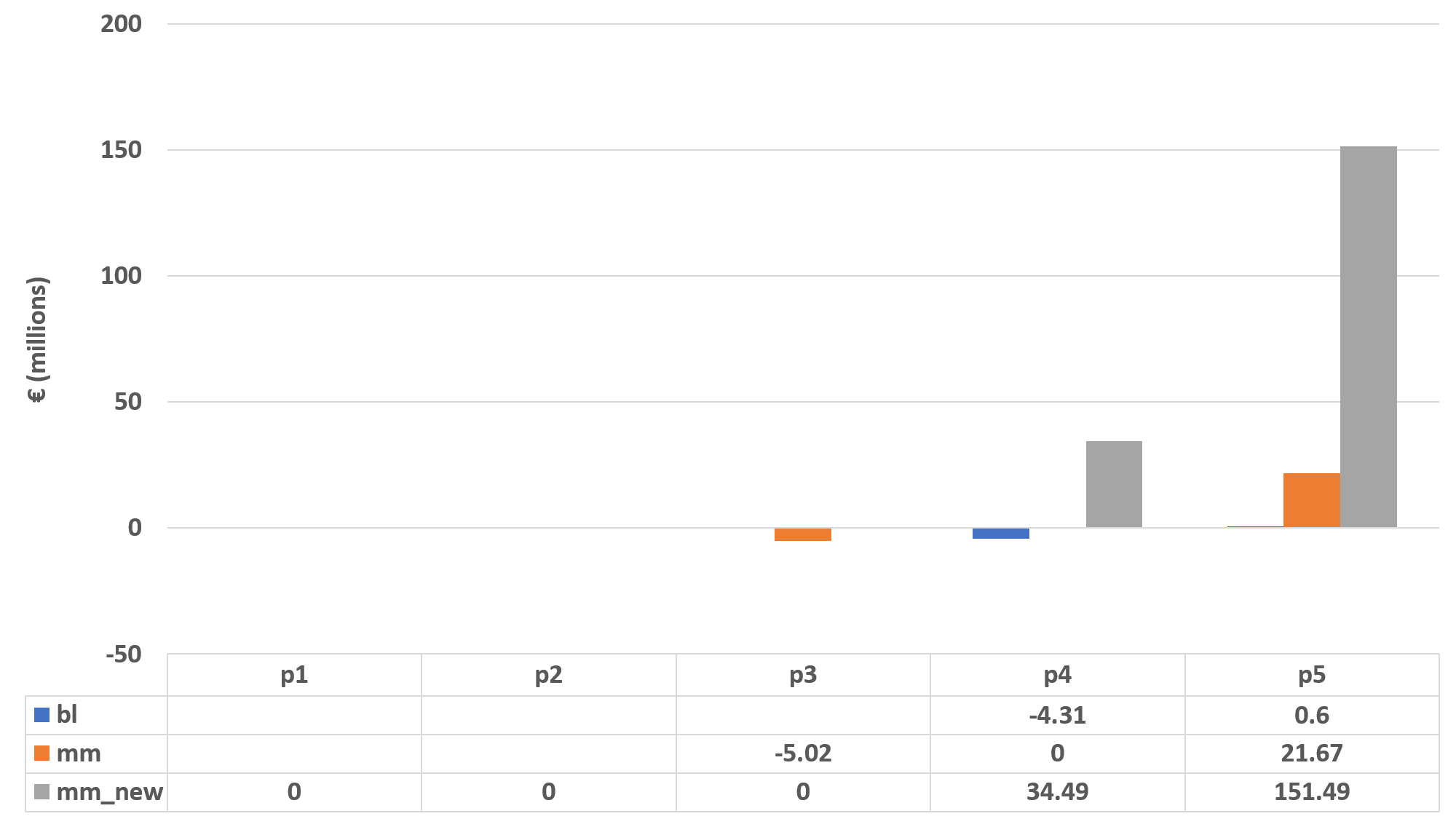}
		\caption{Revenue earned by firm $l=1$ for the first successful iteration.}\label{fig:profit_loss_l1}
\end{figure}

\begin{figure}[htbp]
	\centering
		\includegraphics[width=1\textwidth]{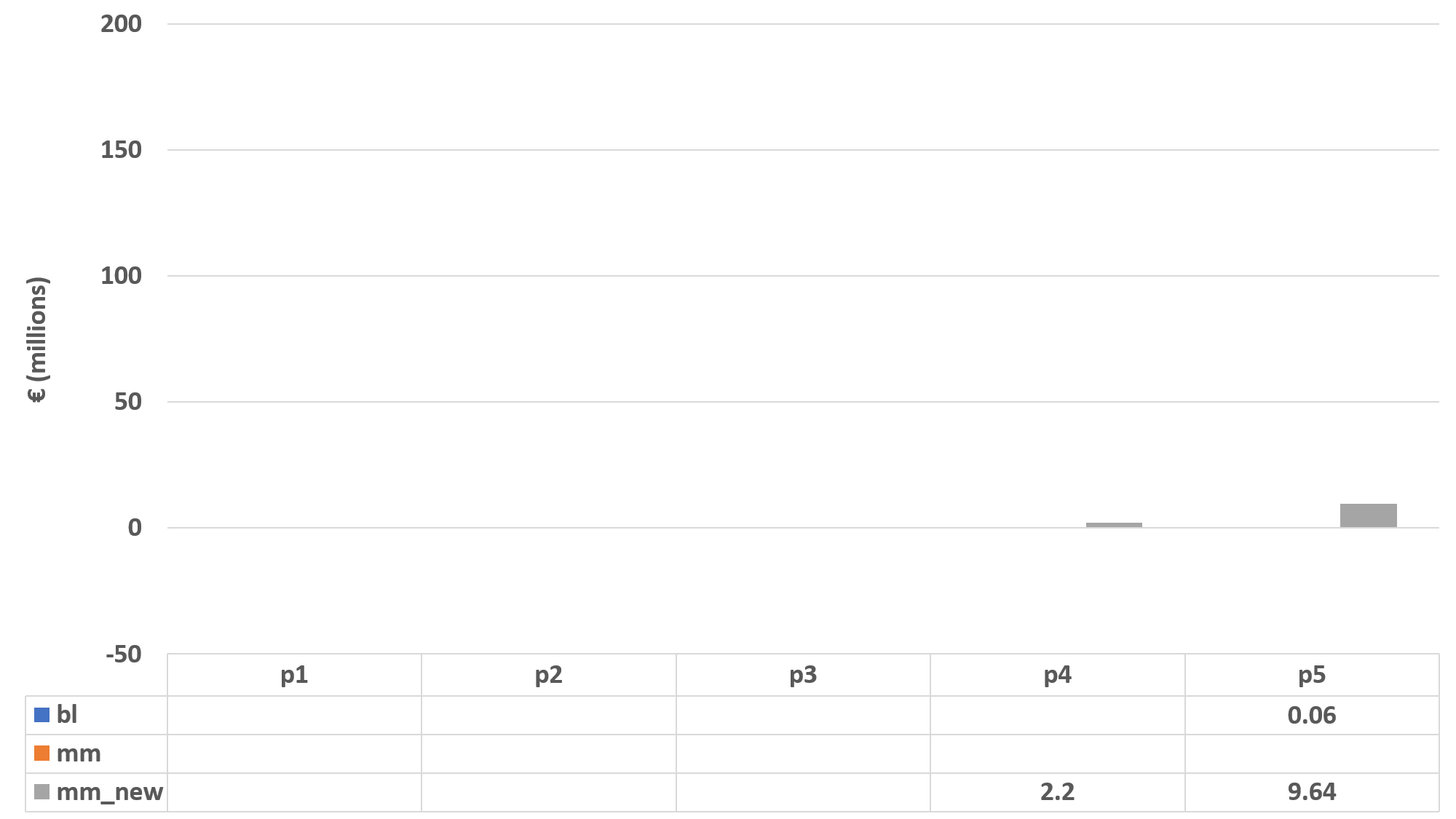}
		\caption{Revenue earned by firm $l=2$ for the first successful iteration.}\label{fig:profit_loss_l2}
\end{figure}

The forward prices in Figure \ref{fig:pices2} can be explained by Figures \ref{fig:gen_SI1} -- \ref{fig:profit_loss_l2}. Figure \ref{fig:gen_SI1} shows the generation mix for the first successful iteration while Figures \ref{fig:profit_loss_l1} and \ref{fig:profit_loss_l2} display the revenues earned/lost by firms $l=1$ and $l=2$, respectively, in each time period for the same iteration. At this equilibrium point, forward prices converged to first time series in Figure \ref{fig:pices2} and firm $l=1$ and firm $l=2$ invested 2765MW and 176MW into new mid-merit generation, respectively. In time periods $p=1$ and $p=2$, only firm $l=1$ and $l=2$'s new mid-merit units are generating leading to forward prices of $\gamma_{p=1}=\gamma_{p=2}=34$, the marginal cost of new mid-merit. Consequently, neither price-making firm earns, nor loses, revenue at these time periods. 

The forward price is the same for $p=3$ but, because the demand curve intercept is higher (see Table \ref{tab:demand}), more generation is needed to meet demand. The increased demand is primarily met by firm $l=2$'s new mid-merit unit. In addition however, firm $l=1$'s existing mid-merit unit generates 403MW. This is despite existing mid-merit having a marginal cost of 41.1. Thus, as Figure \ref{fig:profit_loss_l1} outlines, firm $l=1$ losses revenue at this timepoint. Firm $l=2$ does not earn revenue, nor does it lose revenue, at $p=3$. 

In time period $p=4$, the forward price is 41.1 which is the marginal cost of an existing mid-merit unit. Consequently, all mid-merit units, for firm $l=1$, $l=2$ and $f=3$ are utilised. In addition, firm $l=1$ also utilises is existing baseload despite the marginal cost of exiting baseload being 48.87. The forward price is $\gamma_{p=5}=65.19$ at timestep $p=5$. Because this price is higher than the marginal cost of exiting baseload, both price-making firms utilise their existing baseload units and make a profit from doing so. The two price-making firms use their generation to set $\gamma_{p=5}=65.19$ and hence maximise their respective profits. This forward price allows the two price-making firms to partially recover the investments cost associated with investing in new mid-merit generation. Because they both do not earn any revenue from new mid-merit in timesteps $p=1,2,3$, the remaining investment costs are recovered in timestep $p=4$ where the market price of $\gamma_{p=4}=41.1$ allows both firm $l=1$ and $l=2$ to earn enough revenue from their new mid-merit units to break even on their investments. If either price-making firm adjusted their generation so as to set a price higher than 41.1 in $p=4$ or higher than 65.19 in $p=5$, then the two price-taking firms would invest in new mid-merit generation also, as it would be a profitable decision. The price-making firms prevent this because investment from the price-taking fringe would erode the substantial revenues they earn in timestep $p=5$.

Similarly, it is profit maximising for firm $l=1$ to generate using its existing mid-merit unit, at below marginal cost, in time period $p=3$. If firm $l=1$ did not do this, the remaining demand would be met by firm $f=3$'s existing mid-merit unit, which would drive up the market price and thus, make investing in a new mid-merit a profitable option for both price-taking firms. Again, such a market outcome, is not profit-maximising for firm $l=1$. Instead, it is optimal for firm $l=1$ to take take the small losses in time period $p=3$ so as to prevent the fringe from eroding its large profits in timestep $p=5$. As Figure \ref{fig:profit_loss_l1} shows, firm $l=1$'s revenues from $p=4$ and $p=5$ far exceed its losses from $p=3$. Additionally, in time periods $p=1,2$, it is optimal for firm $l=1$ to ensure the market price is $\gamma_{p=1}=34$. However, in these timesteps, firm $l=1$ does not need its existing mid-merit unit to maintain the price at this level.

These results are in contrast to the prices observed in the perfect competition case; see Figure \ref{fig:pices2}. In the perfect competition setting, firms only utilise a generating unit if the market price is at or above the marginal cost of that unit. Consequently, the market price is set by the marginal cost of the most expensive unit that is generating. Hence, there is no below marginal cost operation of units in time periods $p=3$ and $p=4$, which leads to higher forward prices compared with the oligopoly with a competitive fringe case. Similarly, in time period $p=5$, the forward price in the perfect competition case is set by the most expensive unit that is generating; existing baselaod. In contrast, in the oligopoly with a competitive fringe case, it is optimal for the price-making firms to adjust their generation to ensure the forward price is higher than the perfect competition case.
 
 Similar qualitative results to those in Figures \ref{fig:gen_SI1} and \ref{fig:profit_loss_l1} can be seen in the rest of the successful iterations. The exact level of revenue earned or lost in each timestep, for both price-making firms, varies in a similar manner to Figures \ref{fig:profits_l1} - \ref{fig:invest}.
 
\ref{sec:app_pm} and equation \eqref{eqn:kkt_firm_gen_MCP} shows that when using a MCP model, it is never optimal for a generator to operate  one of its units at below marginal cost.  In contrast, when the EPEC approach of this work is utilised, equation \eqref{eqn:LM_genL} shows that when $gen^{\text{PM}}_{l,t,p}>0$, $\gamma_{p}$ can be less than $C^{\text{GEN}}_{t}$. This is because of the additional $\alpha^{\text{KKT\_GEN}}_{l, ff, tt,p}$ that is in equation \eqref{eqn:LM_genL} but not in equation \eqref{eqn:kkt_firm_gen_MCP}. Moreover, this further highlights the benefit of the EPEC approach and the limitations of the MCP approach when modelling an oligopoly with a competitive fringe and investment decisions.  

The generation levels of price-taking firms $f=3$ and $f=4$ were similar for all equilibria that converged to the first two equilibrium price series. Both price-taking firms utilised their existing mid-merit and peaking units, respectively, to maximum capacity in time period $p=5$ only as this was the only time period where the price was high enough for them to make earn profits. As the equilibrium forward prices converged to one of only three series, the price-taking firms' profits similarly converged to one of three levels. For equilibria that converged on the first price time series in Figure \ref{fig:pices2}, the profits were \euro 17.1M and \euro 0.7M for firms $f=3$ and $f=4$, respectively, while for equilibria with the second time series, the profits were \euro 22.14M and \euro 3.6M, respectively. At the equilibrium point where the third price time series was observed firm $f=3$ also utilised its existing at time period $p=4$ in addition to $p=5$. At this equilibrium point, firm $f=4$ did not generate any electricity as the price was never high enough for them to so. Consequently, firm $f=4$ made zero profits while firm $f=3$ made a profit of \euro 17.1M.

\begin{figure}[htbp]
	\centering
		\includegraphics[width=1\textwidth]{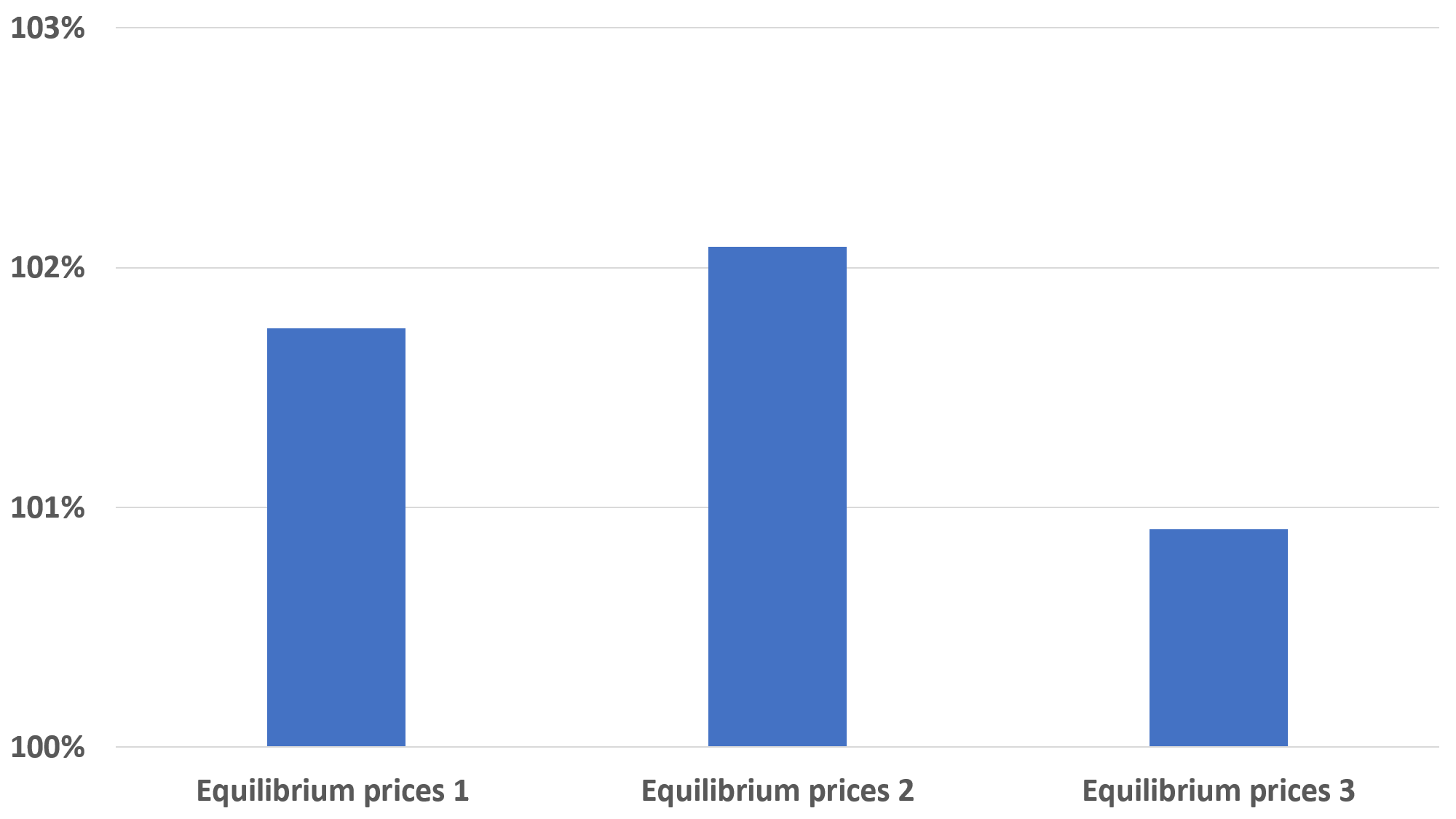}
			\caption{Consumer costs as \% of perfect competition case.}\label{fig:con_costs}
\end{figure}

Figure \ref{fig:con_costs} displays the consumer costs, as defined by the following equation:
\begin{equation}\label{eqn:con_costs}
\sum_{p} \bigg( W_{p}\times \gamma_{p} \times \big( \sum_{ll,tt} gen^{\text{PM}}_{ll,tt,p}+\sum_{ff,tt} gen^{\text{PT}}_{ff,tt,p} \big)\bigg).
\end{equation}
As above, because the equilibrium prices landed at one of three price series, consumer costs also converged to one of three levels. This is because the amount of energy consumed has a fixed relationship with market prices; see market clearing condition \eqref{eqn:MCC}. Figure \ref{fig:con_costs} shows how the consumer costs increase by 1.7\%, 2.1\%, 0.9\% for equilibria that converged to first, second and third time series of forward prices, respectively.

In previous works that use similar data, the presence of price-making behaviour was found to lead to a larger increases in consumer costs \citep{devine2018examining}. However, the ability of the price-taking fringe to invest in new generation motivates the price-making firms to reduce forward market prices in some time periods. While the market prices increase again in subsequent time periods, these consumer cost results shows how the presence of a competitive fringe helps mitigate against the negative effects of market power.

\begin{figure}[htbp]
   \begin{subfigure}[b]{0.5\textwidth}
         \centering
         \includegraphics[width=\textwidth]{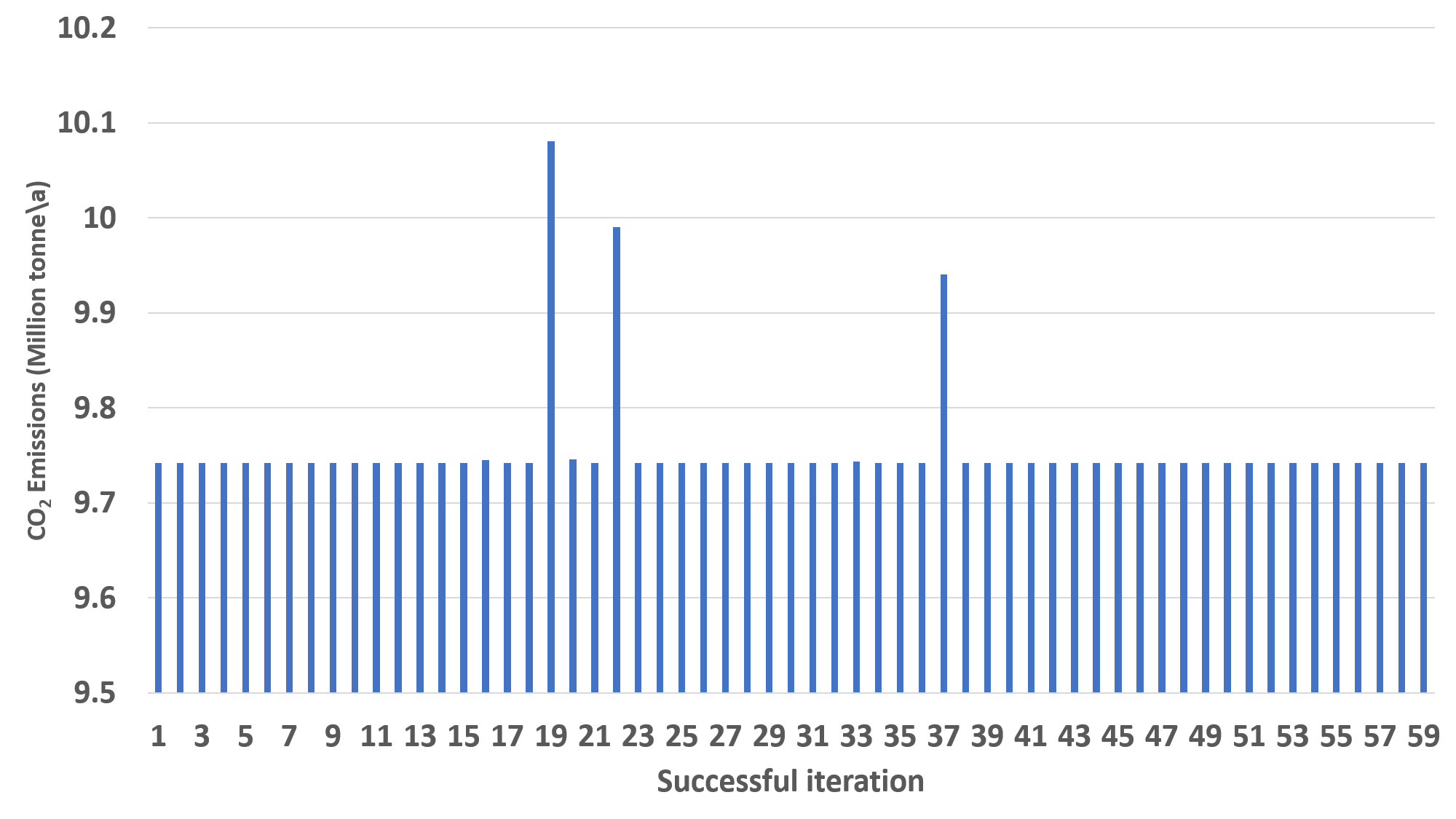}
         \caption{Equilibrium prices 1}
         \label{fig:emissions1}
     \end{subfigure}
      \begin{subfigure}[b]{0.5\textwidth}
         \centering
         \includegraphics[width=\textwidth]{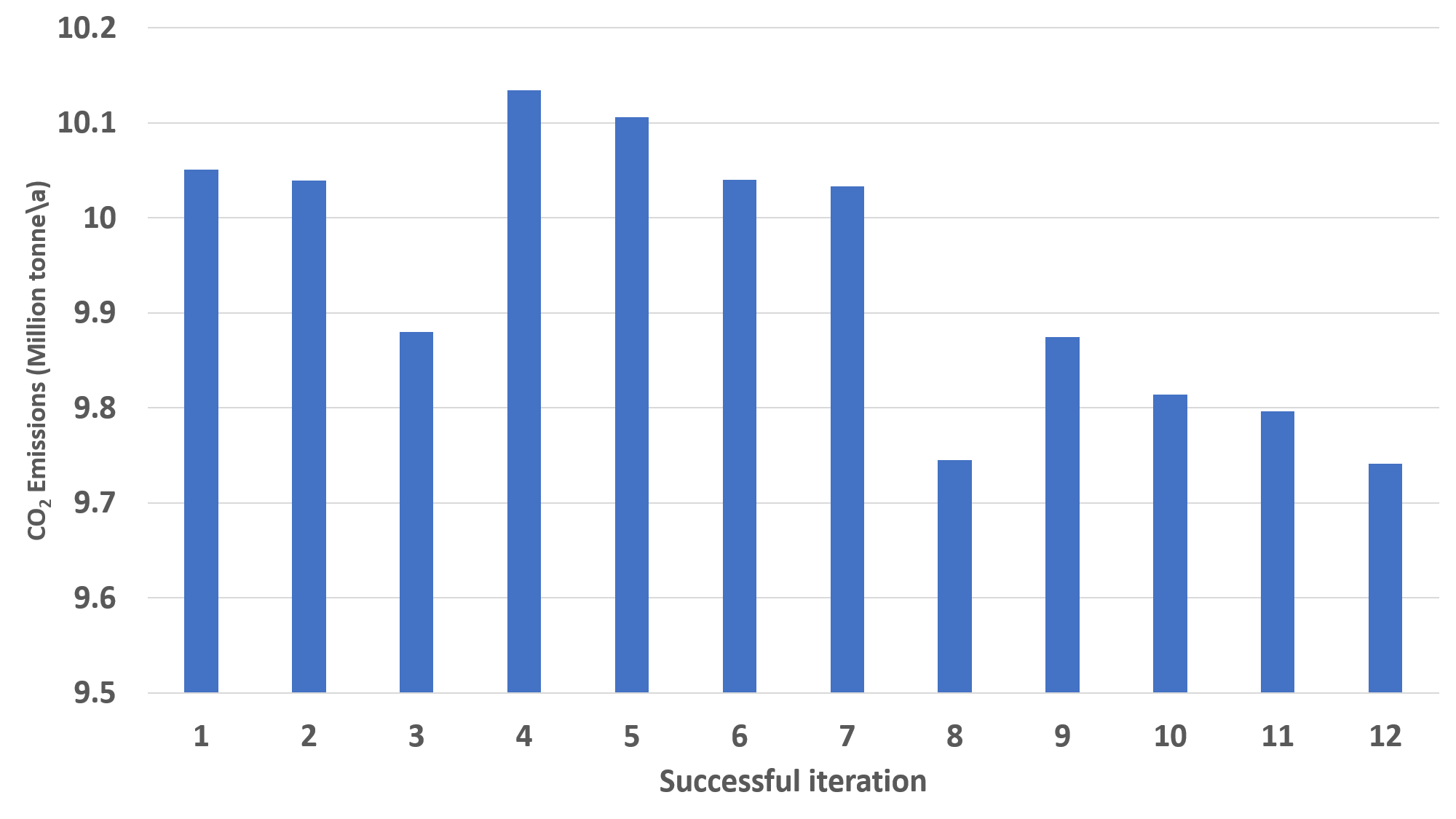}
         \caption{Equilibrium prices 2}
         \label{fig:emissions2}
     \end{subfigure}
     \caption{C0$_2$ emission levels}\label{fig:emissions}
\end{figure}

Figures \ref{fig:emissions1} and \ref{fig:emissions2} display the carbon dioxide emissions level for equilibria that converged at the first and second set of price time series, respectively. These emissions were calculated as follows:
\begin{equation}\label{eqn:emissions}
 \sum_{p,t} \bigg( W_{p}\times E_{t} \times \big( \sum_{ll} gen^{\text{PM}}_{ll,t,p}+\sum_{ff} gen^{\text{PT}}_{ff,p} \big)\bigg),   
\end{equation}
where the parameter $E_{t}$ gives the emissions factor level for technology $t$, as displayed in Table \ref{tab:data_tech}. Figure \ref{fig:emissions} show that, despite equilibrium prices remaining constant across subsets of the equilibria found, the emissions levels varied across the equilibira. This is particularly evident in Figure \ref{fig:emissions2} for equilibria with the second price time series. The reasons behind this results will now be explained.

\subsection{Reasons for multiple equilibria}\label{sec:reasons}
\begin{figure}[htbp]
	\centering
		\includegraphics[width=1\textwidth]{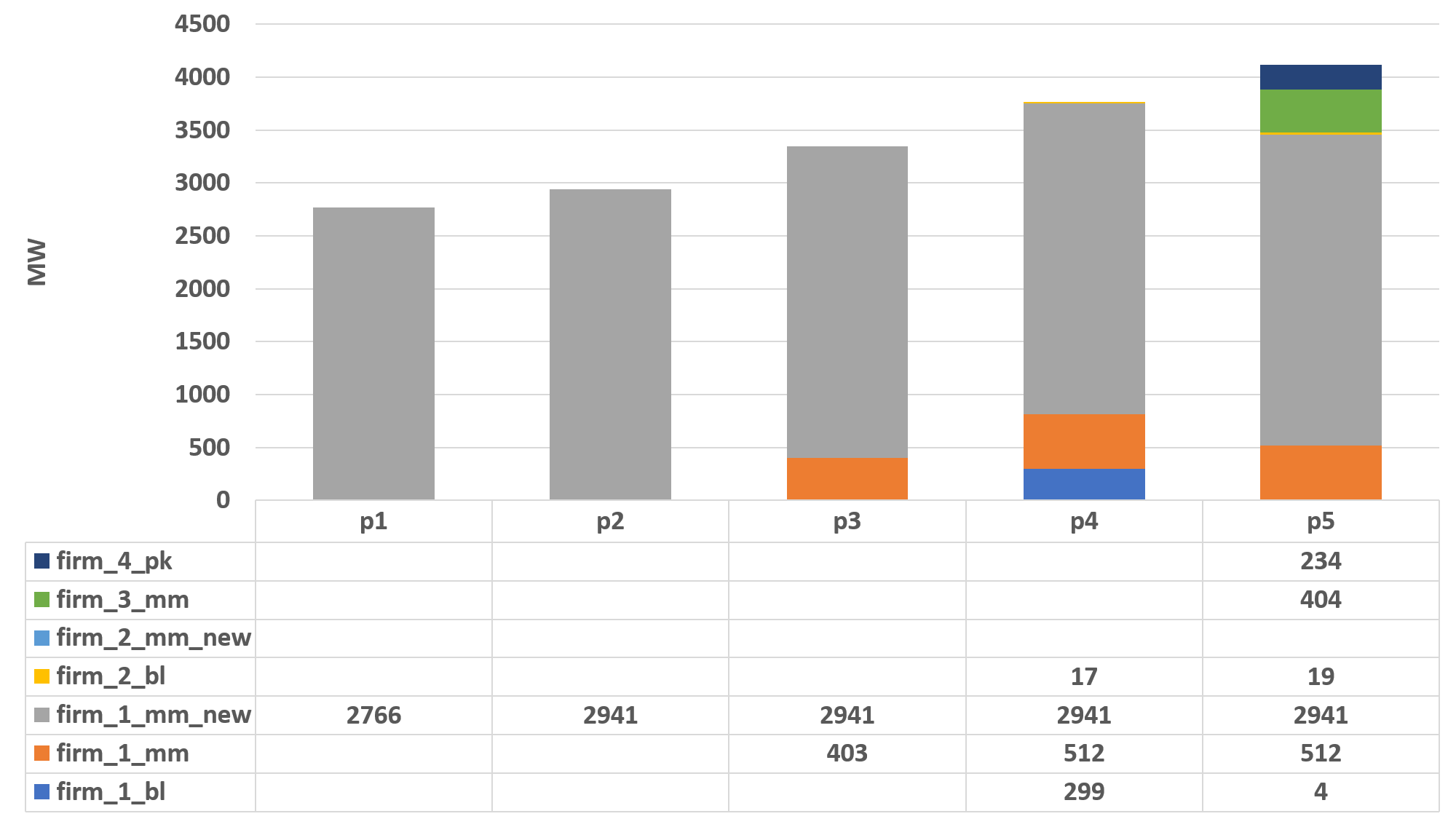}
		\caption{Generation mix for the second successful iteration.}\label{fig:gen_SI2}
\end{figure}

Figures \ref{fig:profits_l1} - \ref{fig:invest} show that there are multiple equilibria to the EPEC presented in this work. In this subsection, we explore the reasons behind this finding. Firstly, we look at two equilibria where the forward prices were the same, i.e., the first equilibrium time series from Figure \ref{fig:pices2}. 

The multiple equilibria are driven by the market's indifference to what firm is providing electricity when firms are generating at the same price. For example, Figures \ref{fig:gen_SI1} and \ref{fig:gen_SI2} display the generation mixes for the first two successful iterations, respectively. In the first successful iteration firms $l=1$ and $l=2$ invest 2765MW and 176MW into new mid-merit generation, respectively. In contrast, in the second successful iteration, they invest 2941MW and 0MW into new mid-merit generation, respectively. 

At time period $p=5$ in the first iteration, firm $l=1$ uses its new and existing mid-merit units to maximum capacity while also generating 21MW from its baseload unit. As the same iteration, firm $l=2$ uses its new mid-merit unit to full capacity and also generates 2MW from its baseload unit. In the second successful iteration, firm $l=1$ decreases its baseload generation at $p=5$ from 21MW to 4MW but increases its generation from new mid-merit from 2765MW to 2941MW. This allows firm $l=1$ to make less profits in Figure \ref{fig:gen_SI2}; firms break even on their new mid-merit investments but make profits from existing baseload generation. Firm $l=2$ increases its baselaod generation from 2MW to 19MW but decreases generation from new mid-merit from 176MW to 0MW. This allows firm $l=1$ to make more profits in Figure \ref{fig:gen_SI2}.

Because the market prices are the same across both equilibria considered, the market is indifferent to whether the electricity comes firm $l=1$'s baseload or mid-merit or from firm $l=2$'s baseload or mid-merit units. Once firm $l$ commits to forward generation decisions, firm $\hat{l} \ne l$ is not willing to adjust its generation levels so as to either increase or decrease forward market price of $\gamma_{p=5}=65.19$. If either price-making firm increased any of the forward prices, then the price-taking firms would invest in new mid-merit generation, as explained in the previous subsection. It is also not profit-maximising for firm $l$ to undercut firm $\hat{l}$ $\ne$ $l$ at a price lower than 65.19. To do so, would mean firm $l$ would make a loss on its new mid-merit investment. Furthermore, if firm $l$ adjusted its generation so as to decrease $\gamma_{p=5}$ by \euro 1, then it would only be able to, at most, increase its generation from existing baseload $\frac{1}{B}=0.11$MW (see market clearing condition \eqref{eqn:MCC}). This is because it would continue to be profitable for firm $\hat{l}$ $\ne$ $l$ to utilise its existing baseload at the reduced price. The small increase in generation opportunity would not make up for the decreased revenues resulting from the reduced price. This paragraph explains why in some of the equilbria found, firm $l=1$ makes less profits than in the perfect competition case. 

Similar market in-differences are also observed in time periods $p=2$-$4$ and in the other 57 successful iterations that converge to the same price time series, thus explaining the multiple equilibria displayed in Figures \ref{fig:profits_l1} - \ref{fig:invest}. In some of other equilibria found, both price-making firms generate significant amounts from their baselaod units in $p=5$, thus preventing each other from generating and investing in new mid-merit generation. Consequently, both firms do not make as large a profit as they otherwise could. Such equilibria are also evident in Figures \ref{fig:profits_l1} - \ref{fig:invest}.  This particular result highlights the absence of collusion between the two price-making firms modelled in this work.

\begin{figure}[htbp]
	\centering
		\includegraphics[width=1\textwidth]{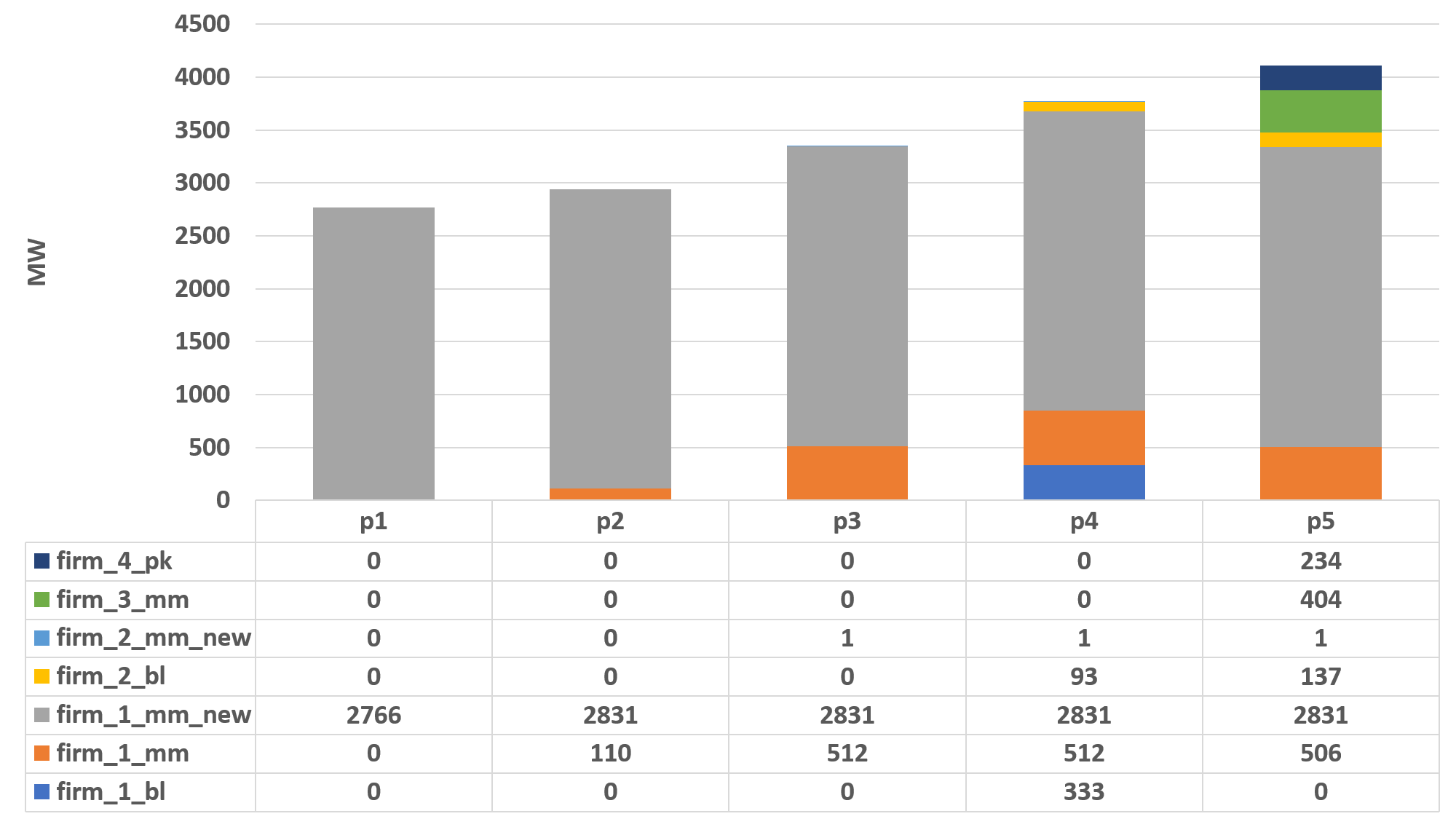}
		\caption{Generation mix for the first successful iteration that results in the second time series for forward prices.}\label{fig:gen_w201}
\end{figure}

\begin{figure}[htbp]
	\centering
		\includegraphics[width=1\textwidth]{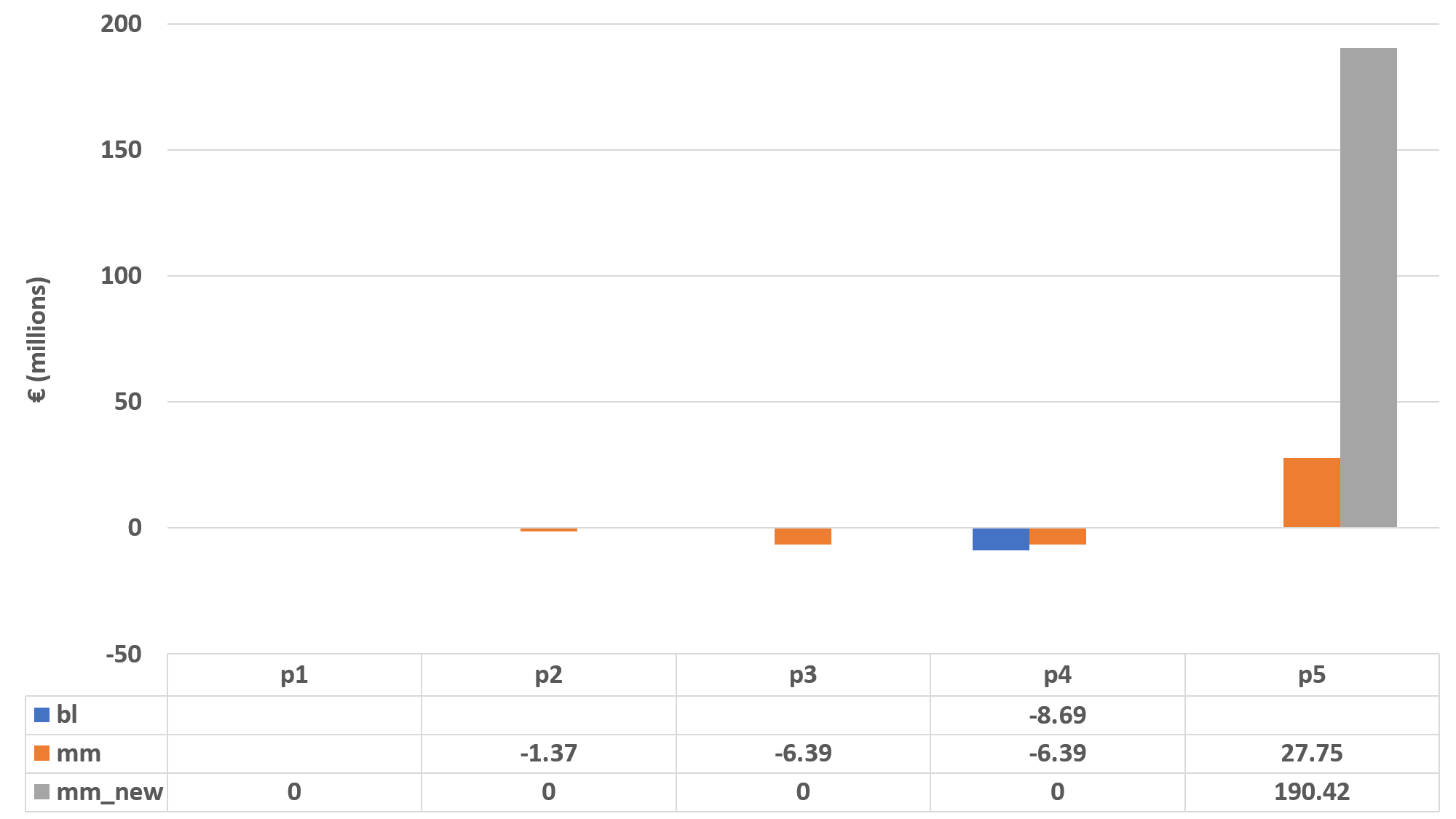}
		\caption{Revenue earned by firm $l=1$ for the first successful iteration that results in the second time series for forward prices.}\label{fig:profit_loss_l1_w201}
\end{figure}

\begin{figure}[htbp]
	\centering
		\includegraphics[width=1\textwidth]{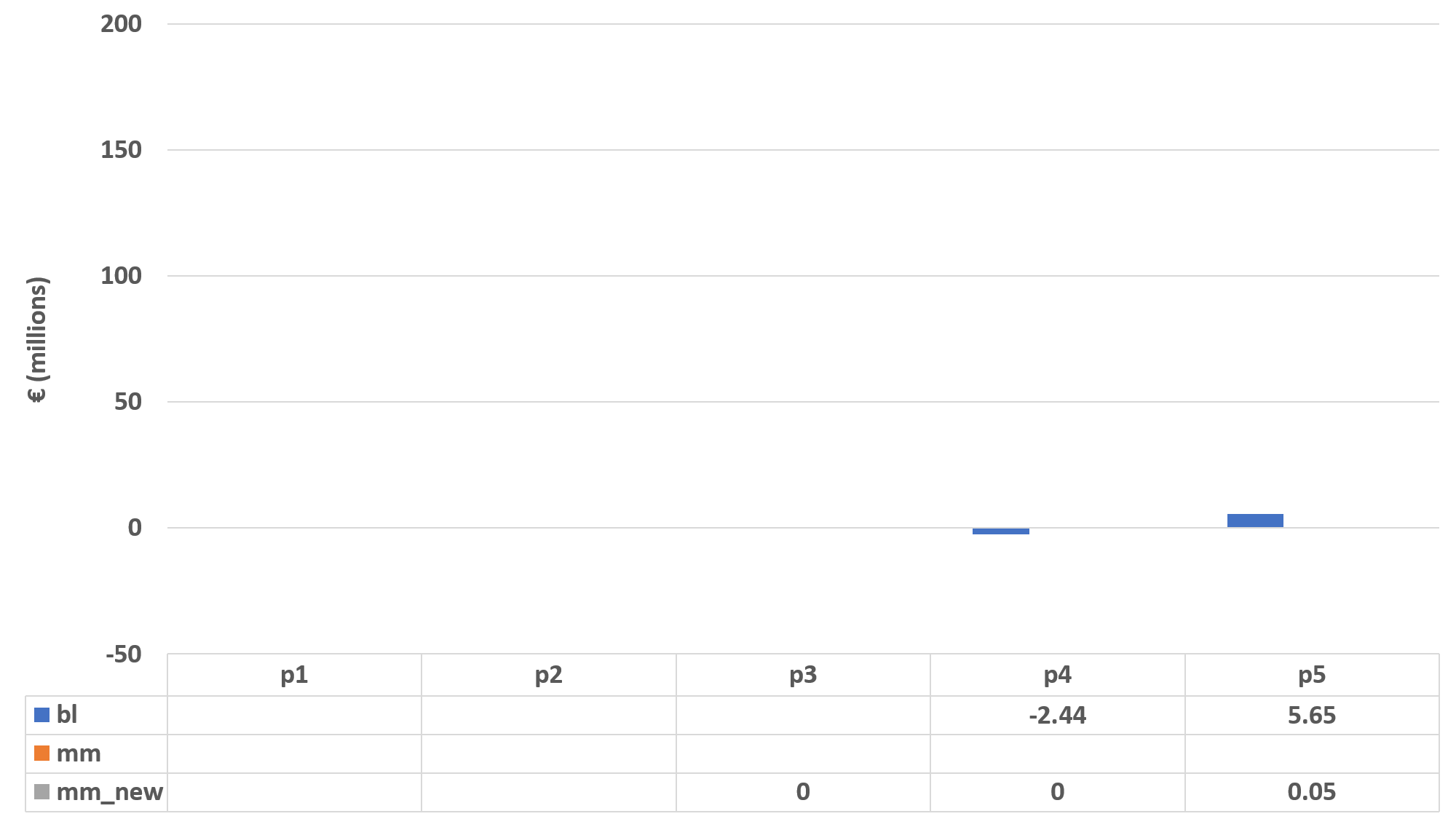}
		\caption{Revenue earned by firm $l=2$ for the first successful iteration that results in the second time series for forward prices.}\label{fig:profit_loss_l2_w201}
\end{figure}

We now examine the differences between two equilibria that converged to different forward price time series. Figure \ref{fig:gen_w201} displays the generation mix for the first successful iteration where the forward prices converged to the second time series in Figure \ref{fig:pices2} while Figures \ref{fig:profit_loss_l1_w201} and \ref{fig:profit_loss_l2_w201} show the revenues for firms $l=1$ and $l=2$, respectively, for the same iteration. Firms $l=1$ and $l=2$ made  profits of \euro 4.92M and \euro 3.21M, respectively, at the equilibrium point. 

In Figure \ref{fig:gen_SI1}, firm $l=2$ generated 17MW from its existing baseload  unit at time period $p=4$. In contrast, in Figure \ref{fig:gen_w201}, firm $l=2$ increased its baseload generation to 93MW at time period $p=4$. Following from market clearing condition \eqref{eqn:MCC}, this lead to the forward price decreasing from $\gamma_{p=4}=41.1$ to $\gamma_{p=4}=31$ between the two time series. This decrease in forward price meant that firm $l=1$ needed to decrease its overall generation in $p=5$ from 3456MW in Figure \ref{fig:gen_SI1} to 3337MW in Figure \ref{fig:gen_w201}. This resulted in a higher forward price in $p=5$ for in the second time series and thus allowed both price-making firms to recover its investment capital costs, despite the decreased price in $p=4$.

Firm $l=2$ cannot make a profit from its existing baseload unit in time period $p=4$ at either $\gamma_{p=4}=41.1$ or $\gamma_{p=4}=31$. Consequently, firm $l=2$ prefers a higher forward price in $p=5$ as this allows it to maximise its profits on its existing baselaod unit; see Figure \ref{fig:profit_loss_l2} compared with Figure \ref{fig:profit_loss_l2_w201}. In contrast, firm $l=1$ prefers the first forward price time series, i.e., a higher price in $p=4$ and a slightly lower price in $p=5$. In time period $p=4$, firm $l=1$ can earn positive revenues from its new mid-merit unit and not make a loss from its existing mid-merit unit in $p=4$, if the forward price is 41.1. In contrast however, firm $l=2$ does not own an exiting mid-merit unit and, in Figures \ref{fig:gen_w201} -- \ref{fig:profit_loss_l2_w201}, only invests in 1MW of new mid-merit generation. Consequently, firm $l=2$ prefers a lower forward price in $p=4$ and a higher price in $p=5$ as this allows firm $l=2$ to maximise its profits from its existing baselaod unit. For both forward price time series, the forward price is not high enough for existing baseload units to earn positive revenues in $p=4$.

In general, the equilibria resulting from the second time series represent equilibria where firm $l=2$ has invested in a relatively small amount of new mid-merit generation, if any at all, but where firm $l=2$ has also made forward generation decisions before firm $l=1$. When firm $l=2$ commits to a large amount of generation in $p=4$, firm $l=1$ must reduce its generation in $p=5$ in order to allow the forward price increase and hence break even on its and firm $l=2$'s new mid-merit investments. 

Interestingly, there is one equilibrium point where firm $l=1$ does not invest in any new technology and consequently, commits to a large amount of generation in $p=4$. This leads to the second equilibrium price time series from Figure \ref{fig:pices2}. This also forces firm $l=2$ to reduce its generation in $p=5$ and motivates it to not make any investment decisions either. As a result, this is the only equilibrium point where the followers make investment decisions; firm $f=3$ invests 2766MW into a new mid-merit facility while firm $f=4$ invest 86MW into the same technology. Because of the generation commitments of the price-making firms set the equilibrium prices, both price-taking firms break exactly even on these investments. Thus, in the model, the price-making firms are indifferent to whether they do the investment at this equilibrium point or the price-taking firms do. However, this indifference may not reflect reality. In the real-world, price-making firms may fear losing their price-making ability if they allow the competitive fringe to invest.  The EPEC model presented in this work does not account for this as the price-making/price-taking characteristics of all firms remain unchanged throughout the model.

Finally, as mentioned above, there was one equilibrium point found where the prices converged to third equilibrium time series in Figure \ref{fig:pices2}. In comparison with the second equilibrium time series, firm $l=2$ commits to investing in new mid-merit generation before $l=1$. However, at this equilibrium point, firm $l=1$ also commits to a large amount of generation in $p=5$ which leads to a reduced price of $\gamma_{p=5}=47.79$, from which its existing mid-merit unit profits from. Consequently, in order for firm $l=2$ to break even on its new mid-merit investment, firm $l=2$ is forced to ensure its generation in $p=4$ is low enough to allow $\gamma_{p=4}=58.59$.

\section{Discussion}\label{sec:dis}
The following summarises the five main findings of our research. Firstly, an Equilibrium Problem with Equilibrium Constraints (EPEC) is a prudent model choice  when modelling an oligopoly with competitive fringe and investments. As outlined in Section \ref{sec:results_mcp}, when investment decisions are included in the model, using a Mixed Complementarity Problem (MCP) can lead to myopic model behaviour and thus contradictory results. Our analysis shows that an EPEC model can overcome this issue and moreover does not require the limiting assumption of conjectural variations. 

Secondly, our results show that may it be optimal for generating firms with market power to occasionally operate some of their generating units at a loss. The driving factor behind this model outcome is the fact we allow both price-making and price-taking firms to make investment decisions. The ability of price-taking firms to invest further into the market motivates the price-making firms to depress prices in some timepoints. This reduces the revenues price-taking firms could make from new investments and thus prevents them from making such investments. Such behaviour would not be captured by MCP or cost-minimisation unit commitment models. Consequently, this result again highlights the suitability of the EPEC modelling approach and the importance of including investment decisions in models of oligopolies with competitive fringes.

Thirdly, the analysis in Section \ref{sec:results_epec} found multiple market equilibria. This led to varied investment decisions, and thus profits, for the price-making firms. These results will be of interest to generating firms, particularly those with market power. Figures \ref{fig:profits_l1} - \ref{fig:invest} highlight the benefit of making investment decisions before other competing price-making firms do so. In fact, our results indicate that if firms do not expand their generation portfolios, then they may face profits lower than they would if the market was perfectly competitive. The multiple equilbira result also indicates that that generation from existing baseload generation may be higher in some equilibria compared with others. Such market outcomes will be of interest to energy policymakers who are concerned about carbon emission levels. Older baseload generators tend to be coal-based and thus emit higher levels of carbon. Consequently, while the market may be indifferent to where the electricity comes from, policymakers may seek to put measures in places to encourage the equilibrium outcomes where existing baseload generation is reduced. 
 
Fourthly, Figure \ref{fig:con_costs} showed that the presence of market power increases consumer costs by 1\% -- 2\%. While this outcome is not surprising, the level is relatively small compared with the literature. For instance, using similar data, \cite{devine2019role} estimate market power in an oligopoly with competitive fringe context could double consumer costs compared with a perfectly competitive market. However, \cite{devine2019role} do not include investment decisions in their model. Thus, this result again highlights the impact of including investment decisions in models of oligopolies with competitive fringes. It also highlights the importance for policymakers to encourage new entrants into electricity markets and, moreover, the benefits of encouraging smaller generating firms to expand their portfolios, or at least threaten to.

Finally, as the literature details \citep{pozo2017basic}, solving EPEC problems can be computational challenging. In this work we utilised the method outlined in \cite{leyffer2010solving} to obtain an initial starting point solution to our algorithm. Using this approach our algorithm successfully found an equilibrium from 72 of the 200 iterations attempted. When instead we used a random initial starting point solution, we found an equilibrium from only 2 of the 2000 iterations attempted. 

Critically reflecting on our approach, we wish to acknowledge some limitations. Firstly, because EPEC problems are challenging to solve, we choose the relatively small number of five timesteps. These represented hours in summer low demand, summer high demand, winter low demand, winter high demand and winter peak demand. Thus, the net demand intercept values represent average values for these timesteps. In reality, particularly in systems with a large amount of renewables, these intercept values will fluctuate from hour to hour. As a result, the average values may over- or under-estimate the total profits each generating firm could make in each time period. This would impact investments decisions and consumer costs.

Secondly, we did not account for any stochasticity in the model. Due to the intermittent and uncertain nature of wind energy, stochasticity is a feature of many electricity market models. Such stochasticity is typically introduced by making generation capacity scenario-dependent \citep{lynch2019role}. Deterministic capacity values may also over- or under-estimate the profits each firm may make in each timestep. However, while including further timesteps and stochastic capacity values would most likely affect the exact numbers presented in Section \ref{sec:results_epec}, we do not anticipate them changing the qualitative findings discussed above. 

Finally, we did not consider a capacity market as part of the market modelled in this work. Capacity payments exists when firms get paid for simply owning generation units and making them available to the grid. Capacity payments do not depend on the extent that the unit(s) are utilised. Regulators and policymakers include such payments so as to ensure security of supply \citep{lynch2017investment}. The market we considered was an 'energy-only' market, where the generating firms only get paid on the basis if how much they generate. A capacity market can affect the level of investment into new generation. Future research activities will address each of these modelling limitations.

\section{Conclusion}\label{sec:con}
In this paper, we developed a novel mathematical model of an imperfect electricity market, one that is characterised by an oligopoly with a competitive fringe. We modelled two types of generating firms; price-making firms, who have market power, and price-taking firms who do not. All firms had both investment and forward generation decisions. The model took the form of an Equilibrium Problem with Equilibrium Constraints (EPEC), which finds an equilibrium of multiple bi-level optimisation problems. The bi-level formulation allowed the optimisation problems of the price-taking firms to be embedded into the optimisation problems of the price-making firms. This enabled the price-making firms to correctly anticipate the optimal reactions of the price-taking firms to their decisions. We applied the model to data representative of the Irish power system for 2025. 

To solve the EPEC problem, we utilised the Gauss-Seidel algorithm. Furthermore, we found the computational efficiency of the algorithm was improved when the algorithm's starting point solution was provided by the approach detailed in \cite{leyffer2010solving}. Overall, we found that an EPEC problem is a prudent model choice  when modelling an oligopoly with competitive fringe. This is because it overcomes modelling issues previously found in the literature and requires a fewer limiting assumptions. 

The model found multiple equilibria. This was due to the market's indifference to which price-making firm generates electricity. Although consumer costs were found to be relatively constant across the equilibria found, this result is important to energy policymakers who who wish to avoid equilibrium outcomes that higher carbon emission levels. 

We also observed that it may be strategically optimal for price-making firms to occasionally to generate at a price that is lower than their marginal cost. This is because we incorporated investment decisions into the optimisation problems of both types of generating firms. Consequently, the price-making firms seek to depress prices occasionally so as to discourage the fringe from investing further into the market. Furthermore, we found that consumer costs only decreased by 1\% -- 2\% when market power was removed from the model.

In future research, we will study the effects of increasing the number of timesteps in the model. Moreover, we will explore the impact stochasticity, particularly from wind generation, would have. In addition, future research will analyse how the introduction of a capacity market would affect equilibrium outcomes.

\section*{Acknowledgements}

M. T. Devine acknowledges funding from Science Foundation Ireland (SFI) under the SFI Strategic Partnership Programme Grant number SFI/15/SPP/E3125. S. Siddiqui acknowledges funding by NSF Grant \#1745375 [EAGER: SSDIM: Generating Synthetic Data on Interdependent Food, Energy, and Transportation Networks via Stochastic, Bi-level Optimization]. The authors also sincerely thank Dr. M. Lynch and Dr. S. Lyons from the Economic and Social Research Institute (ESRI) in Dublin and who provided invaluable feedback and advice.

\section*{References}
\bibliography{bibtex_epec.bib}

\begin{thebibliography}{44}
\expandafter\ifx\csname natexlab\endcsname\relax\def\natexlab#1{#1}\fi
\expandafter\ifx\csname url\endcsname\relax
  \def\url#1{\texttt{#1}}\fi
\expandafter\ifx\csname urlprefix\endcsname\relax\def\urlprefix{URL }\fi

\bibitem[{Ansari(2017)}]{ansari2017opec}
Ansari, D., 2017. Opec, saudi arabia, and the shale revolution: Insights from
  equilibrium modelling and oil politics. Energy Policy 111, 166--178.

\bibitem[{Baltensperger et~al.(2016)Baltensperger, F{\"u}chslin, Kr{\"u}tli,
  and Lygeros}]{baltensperger2016multiplicity}
Baltensperger, T., F{\"u}chslin, R.~M., Kr{\"u}tli, P., Lygeros, J., 2016.
  Multiplicity of equilibria in conjectural variations models of natural gas
  markets. European Journal of Operational Research 252~(2), 646--656.

\bibitem[{Bertsch et~al.(2018)Bertsch, Devine, Sweeney, and
  Parnell}]{bertsch2018analysing}
Bertsch, V., Devine, M.~T., Sweeney, C., Parnell, A.~C., 2018. Analysing
  long-term interactions between demand response and different electricity
  markets using a stochastic market equilibrium model. Tech. rep., ESRI Working
  Paper.

\bibitem[{Bushnell et~al.(2008)Bushnell, Mansur, and
  Saravia}]{bushnell2008vertical}
Bushnell, J.~B., Mansur, E.~T., Saravia, C., 2008. Vertical arrangements,
  market structure, and competition: An analysis of restructured us electricity
  markets. American Economic Review 98~(1), 237--66.

\bibitem[{Devine and Bertsch(2018)}]{devine2018examining}
Devine, M.~T., Bertsch, V., 2018. Examining the benefits of load shedding
  strategies using a rolling-horizon stochastic mixed complementarity
  equilibrium model. European Journal of Operational Research 267~(2),
  643--658.

\bibitem[{Devine and Bertsch(2019)}]{devine2019role}
Devine, M.~T., Bertsch, V., 2019. The role of demand response in mitigating
  market power - a quantitative analysis using a stochastic market equilibrium
  model. Tech. rep., ESRI Working Paper.

\bibitem[{Devine et~al.(2019)Devine, Nolan, Lynch, and
  O’Malley}]{devine2019effect}
Devine, M.~T., Nolan, S., Lynch, M.~{\'A}., O’Malley, M., 2019. The effect of
  demand response and wind generation on electricity investment and operation.
  Sustainable Energy, Grids and Networks 17, 100190.

\bibitem[{Di~Cosmo and Hyland(2013)}]{di2013carbon}
Di~Cosmo, V., Hyland, M., 2013. Carbon tax scenarios and their effects on the
  irish energy sector. Energy Policy 59, 404--414.

\bibitem[{Egging(2013)}]{egging2013benders}
Egging, R., 2013. Benders decomposition for multi-stage stochastic mixed
  complementarity problems--applied to a global natural gas market model.
  European Journal of Operational Research 226~(2), 341--353.

\bibitem[{Egging et~al.(2008)Egging, Gabriel, Holz, and
  Zhuang}]{egging2008complementarity}
Egging, R., Gabriel, S.~A., Holz, F., Zhuang, J., 2008. A complementarity model
  for the european natural gas market. Energy policy 36~(7), 2385--2414.

\bibitem[{Egging and Holz(2016)}]{egging2016risks}
Egging, R., Holz, F., 2016. Risks in global natural gas markets: investment,
  hedging and trade. Energy Policy 94, 468--479.

\bibitem[{EirGrid(2016)}]{eirgrid2016}
EirGrid, 2016. All-island generation capacity statement 2016-2025.

\bibitem[{Facchinei and Pang(2007)}]{facchinei2007finite}
Facchinei, F., Pang, J.-S., 2007. Finite-dimensional variational inequalities
  and complementarity problems. Springer Science \& Business Media.

\bibitem[{Fanzeres et~al.(2019)Fanzeres, Ahmed, and
  Street}]{fanzeres2019robust}
Fanzeres, B., Ahmed, S., Street, A., 2019. Robust strategic bidding in
  auction-based markets. European Journal of Operational Research 272~(3),
  1158--1172.

\bibitem[{Fortuny-Amat and McCarl(1981)}]{fortuny1981representation}
Fortuny-Amat, J., McCarl, B., 1981. A representation and economic
  interpretation of a two-level programming problem. Journal of the operational
  Research Society 32~(9), 783--792.

\bibitem[{Gabriel et~al.(2012)Gabriel, Conejo, Fuller, Hobbs, and
  Ruiz}]{gabriel2012complementarity}
Gabriel, S.~A., Conejo, A.~J., Fuller, J.~D., Hobbs, B.~F., Ruiz, C., 2012.
  Complementarity modeling in energy markets. Vol. 180. Springer Science \&
  Business Media.

\bibitem[{Gabriel et~al.(2009)Gabriel, Zhuang, and Egging}]{gabriel2009solving}
Gabriel, S.~A., Zhuang, J., Egging, R., 2009. Solving stochastic
  complementarity problems in energy market modeling using scenario reduction.
  European Journal of Operational Research 197~(3), 1028--1040.

\bibitem[{Guo et~al.(2019)Guo, Chen, Xia, and Kang}]{guo2019electricity}
Guo, H., Chen, Q., Xia, Q., Kang, C., 2019. Electricity wholesale market
  equilibrium analysis integrating individual risk-averse features of
  generation companies. Applied Energy 252, 113443.

\bibitem[{Habibian et~al.(2019)Habibian, Downward, and
  Zakeri}]{habibian2019multistage}
Habibian, M., Downward, A., Zakeri, G., 2019. Multistage stochastic demand-side
  management for price-making major consumers of electricity in a co-optimized
  energy and reserve market. European Journal of Operational Research.

\bibitem[{Haftendorn and Holz(2010)}]{haftendorn2010modeling}
Haftendorn, C., Holz, F., 2010. Modeling and analysis of the international
  steam coal trade. The Energy Journal, 205--229.

\bibitem[{Hu and Ralph(2007)}]{hu2007using}
Hu, X., Ralph, D., 2007. Using epecs to model bilevel games in restructured
  electricity markets with locational prices. Operations research 55~(5),
  809--827.

\bibitem[{Huppmann(2013)}]{huppmann2013endogenous}
Huppmann, D., 2013. Endogenous shifts in opec market power: a stackelberg
  oligopoly with fringe.

\bibitem[{Huppmann and Egerer(2015)}]{huppmann2015national}
Huppmann, D., Egerer, J., 2015. National-strategic investment in european power
  transmission capacity. European Journal of Operational Research 247~(1),
  191--203.

\bibitem[{Huppmann and Egging(2014)}]{huppmann2014market}
Huppmann, D., Egging, R., 2014. Market power, fuel substitution and
  infrastructure--a large-scale equilibrium model of global energy markets.
  Energy 75, 483--500.

\bibitem[{Jin and Ryan(2013)}]{jin2013tri}
Jin, S., Ryan, S.~M., 2013. A tri-level model of centralized transmission and
  decentralized generation expansion planning for an electricity market—part
  i. IEEE Transactions on Power Systems 29~(1), 132--141.

\bibitem[{Kazempour et~al.(2013)Kazempour, Conejo, and
  Ruiz}]{kazempour2013generation}
Kazempour, S.~J., Conejo, A.~J., Ruiz, C., 2013. Generation investment
  equilibria with strategic producers—part i: Formulation. IEEE Transactions
  on Power Systems 28~(3), 2613--2622.

\bibitem[{Kazempour and Zareipour(2013)}]{kazempour2013equilibria}
Kazempour, S.~J., Zareipour, H., 2013. Equilibria in an oligopolistic market
  with wind power production. IEEE Transactions on Power Systems 29~(2),
  686--697.

\bibitem[{Lee(2016)}]{lee2016nash}
Lee, C.-Y., 2016. Nash-profit efficiency: A measure of changes in market
  structures. European Journal of Operational Research 255~(2), 659--663.

\bibitem[{Leyffer and Munson(2010)}]{leyffer2010solving}
Leyffer, S., Munson, T., 2010. Solving multi-leader--common-follower games.
  Optimisation Methods \& Software 25~(4), 601--623.

\bibitem[{Lynch et~al.(2019{\natexlab{a}})Lynch, Devine, and
  Bertsch}]{lynch2019role}
Lynch, M., Devine, M.~T., Bertsch, V., 2019{\natexlab{a}}. The role of
  power-to-gas in the future energy system: Market and portfolio effects.
  Energy 185, 1197--1209.

\bibitem[{Lynch and Devine(2017)}]{lynch2017investment}
Lynch, M.~A., Devine, M.~T., 2017. Investment vs. refurbishment: examining
  capacity payment mechanisms using stochastic mixed complementarity problems.
  The Energy Journal 38~(2).

\bibitem[{Lynch et~al.(2019{\natexlab{b}})Lynch, Nolan, Devine, and
  O’Malley}]{lynch2019impacts}
Lynch, M.~{\'A}., Nolan, S., Devine, M.~T., O’Malley, M., 2019{\natexlab{b}}.
  The impacts of demand response participation in capacity markets. Applied
  Energy 250, 444--451.

\bibitem[{Moiseeva et~al.(2017)Moiseeva, Wogrin, and
  Hesamzadeh}]{moiseeva2017generation}
Moiseeva, E., Wogrin, S., Hesamzadeh, M.~R., 2017. Generation flexibility in
  ramp rates: Strategic behavior and lessons for electricity market design.
  European Journal of Operational Research 261~(2), 755--771.

\bibitem[{Pozo and Contreras(2011)}]{pozo2011finding}
Pozo, D., Contreras, J., 2011. Finding multiple nash equilibria in pool-based
  markets: A stochastic epec approach. IEEE Transactions on Power Systems
  26~(3), 1744--1752.

\bibitem[{Pozo et~al.(2013)Pozo, Contreras, and Sauma}]{pozo2013if}
Pozo, D., Contreras, J., Sauma, E., 2013. If you build it, he will come:
  Anticipative power transmission planning. Energy Economics 36, 135--146.

\bibitem[{Pozo et~al.(2017)Pozo, Sauma, and Contreras}]{pozo2017basic}
Pozo, D., Sauma, E., Contreras, J., 2017. Basic theoretical foundations and
  insights on bilevel models and their applications to power systems. Annals of
  Operations Research 254~(1-2), 303--334.

\bibitem[{Pozo et~al.(2012)Pozo, Sauma, and Contreras}]{pozo2012three}
Pozo, D., Sauma, E.~E., Contreras, J., 2012. A three-level static milp model
  for generation and transmission expansion planning. IEEE Transactions on
  Power systems 28~(1), 202--210.

\bibitem[{Ruiz et~al.(2011)Ruiz, Conejo, and Smeers}]{ruiz2011equilibria}
Ruiz, C., Conejo, A.~J., Smeers, Y., 2011. Equilibria in an oligopolistic
  electricity pool with stepwise offer curves. IEEE Transactions on Power
  Systems 27~(2), 752--761.

\bibitem[{Steeger and Rebennack(2017)}]{steeger2017dynamic}
Steeger, G., Rebennack, S., 2017. Dynamic convexification within nested benders
  decomposition using lagrangian relaxation: An application to the strategic
  bidding problem. European Journal of Operational Research 257~(2), 669--686.

\bibitem[{Walsh et~al.(2016)Walsh, Malaguzzi~Valeri, and
  Di~Cosmo}]{walsh2016strategic}
Walsh, D., Malaguzzi~Valeri, L., Di~Cosmo, V., 2016. Strategic bidding, wind
  ownership and regulation in a decentralised electricity market.

\bibitem[{Wogrin et~al.(2012)Wogrin, Barqu{\'\i}n, and
  Centeno}]{wogrin2012capacity}
Wogrin, S., Barqu{\'\i}n, J., Centeno, E., 2012. Capacity expansion equilibria
  in liberalized electricity markets: an epec approach. IEEE Transactions on
  Power Systems 28~(2), 1531--1539.

\bibitem[{Wogrin et~al.(2013)Wogrin, Centeno, and
  Barquin}]{wogrin2013generation}
Wogrin, S., Centeno, E., Barquin, J., 2013. Generation capacity expansion
  analysis: Open loop approximation of closed loop equilibria. IEEE
  Transactions on Power Systems 28~(3), 3362--3371.

\bibitem[{Ye et~al.(2017)Ye, Papadaskalopoulos, and
  Strbac}]{ye2017investigating}
Ye, Y., Papadaskalopoulos, D., Strbac, G., 2017. Investigating the ability of
  demand shifting to mitigate electricity producers’ market power. IEEE
  Transactions on Power Systems 33~(4), 3800--3811.

\bibitem[{Zerrahn and Huppmann(2017)}]{zerrahn2017network}
Zerrahn, A., Huppmann, D., 2017. Network expansion to mitigate market power.
  Networks and Spatial Economics 17~(2), 611--644.

\end{thebibliography}

\appendix
\section{Alternative Price-making firm $l$'s problem}\label{sec:app_pm}
When the problem is solved as a Mixed Complementarity Problem (MCP), price-making firm $l$'s optimisation problem takes the following form, where all variables an parameters are as defined previously:
\begin{equation}\label{eqn:Lgen_obj1_MCP}
\begin{split}
\max_{\substack{ 
gen^{\text{PM}}_{l,t,p}, inv^{\text{PM}}_{l,t}\\
gen^{\text{PT}}_{f,t,p}, inv^{\text{PT}}_{f,t}\\
\gamma_{p}, \lambda^{\text{PT}}_{f, t,p}
 }}\>\>
&\sum_{t,p} W_{p}\times gen^{\text{PM}}_{l, t,p} \times \big (  \gamma_{p}   -  C^{\text{GEN}}_{t}\big )-\sum_{t} IC^{\text{GEN}}_{t}\times inv^{\text{PM}}_{l,t}.
\end{split}
\end{equation}
subject to:
\begin{eqnarray}
 gen^{\text{PM}}_{l, t,p} &\leq& CAP^{\text{PM}}_{l,t}+inv^{\text{PM}}_{l,t}, \>\> \forall t,p. \label{eqn:Lgen_energy_con_MCP}
 \end{eqnarray}
The KKT conditions associated with this optimisation problem are
\begin{eqnarray}
 0\leq gen^{\text{PM}}_{l, t,p} &\perp& -W_{p}\times\big(\gamma_{p}  + \frac{\partial \gamma_{p}}{\partial  gen^{\text{PM}}_{l, t,p}} \times gen^{\text{PM}}_{l, t,p}
-C^{\text{GEN}}_{t}\big)+\lambda^{\text{PM}}_{l, t,p} \geq 0, \>\> \forall t,p,\label{eqn:kkt_firm_gen_MCP}\\
0 \leq inv^{\text{PM}}_{l, t} &\perp& IC^{\text{GEN}}_{t}-\sum_{p}\lambda^{\text{PM}}_{l, t,p} \geq 0, \>\> \forall t,\label{eqn:kkt_firm_inv_MCP}\\
0 \leq \lambda^{\text{PM},1}_{l,t,p} &\perp& -gen^{\text{PM}}_{l, t,p}+CAP^{\text{PM}}_{l,t}+inv^{\text{PM}}_{l,t} \geq 0, \>\> \forall t,p,\label{eqn:kkt_firm_lambda_MCP}
\end{eqnarray} 
where 
\begin{equation}\label{eqn:delta_gamma}
    \frac{\partial \gamma_{p}}{\partial  gen^{\text{PM}}_{l, t,p}}=-CV_{l} \times B ,\>\> \forall l,t,p,
\end{equation}
is determined via market clearing condition \eqref{eqn:MCC}. Furthermore, the parameter $CV_{l} \in [0,1]$ represents the Conjectural Variation associated with firm $l$. When firm $l$'s problem is described by equations \eqref{eqn:Lgen_obj1_MCP} and \eqref{eqn:Lgen_energy_con_MCP}, it is a convex optimisation problem and hence the KKT conditions \eqref{eqn:kkt_firm_gen_MCP} - \eqref{eqn:kkt_firm_lambda_MCP} are both necessary and sufficient for optimality \citep{gabriel2012complementarity}.

When the overall market problem is solved as a MCP, the problem consists of the market clearing condition \eqref{eqn:MCC},  the price-taking firms' KKT conditions (equations \eqref{eqn:kkt_firm_gen} - \eqref{eqn:kkt_firm_lambda}) and the KKT conditions for all price-making firms (equations \eqref{eqn:kkt_firm_gen_MCP} - \eqref{eqn:kkt_firm_lambda_MCP}).

It is important to note that when $gen^{\text{PM}}_{l, t,p}>0$, then condition \eqref{eqn:kkt_firm_gen_MCP} is only satisfied if $\gamma_{p} \geq C^{\text{GEN}}_{t}$. Thus, when the above MCP is used in Section \ref{sec:results_mcp}, not one generating unit will operate at below marginal cost. This is in contrast to the EPEC analysis in Section \ref{sec:results_epec} and highlights a further limitation of the MCP modelling approach.

\end{document}